\documentclass[11pt,leqno,openany]{article}
\usepackage{amsmath,amssymb,amsthm,textcomp}
\usepackage{colortbl}
\usepackage[english]{babel}
\usepackage{dsfont}
\newcommand{\bfm}[1]{\boldsymbol{#1}}

\newcounter{ass}

\newtheorem{theorem}{Theorem}[section]
\newtheorem{corollary}[theorem]{Corollary}
\newtheorem{definition}[theorem]{Definition}
\newtheorem{lemma}[theorem]{Lemma}
\newtheorem{proposition}[theorem]{Proposition}

\newtheorem{example}{Example}

\newtheorem{hypothesis}[ass]{Assumption}
\newtheorem{remark}[theorem]{Remark}

\def \H{ \mathbb H}
\def \M{\mathbb M}
\def \R{ \mathbb R }
\def \reel{ \mathbb R }
\def \nat{ \mathbb N}

\def \one{ {\rm 1}\mkern-4.5mu{\rm I} }

\def \E{ \mathbb  E }

\def \P{ \mathbb P  }

\begin{document}
\title{Scaling limits for symmetric It\^o-L\'evy  processes in random medium} 
\author{R\'emi Rhodes${}^{(a)}$ and Vincent Vargas${}^{(a)}$}
\date{}
\maketitle 
{\footnotesize (a)  Address: Universit{\'e} Paris-Dauphine, CEREMADE, Place du Marchal De Lattre de Tassigny, 75775 Paris Cedex 16, France. 
Phone: (33)(0)1 44 05 48 51. E-mail: {\tt rhodes@ceremade.dauphine.fr}} 
\begin{abstract}
We are concerned with scaling limits of the solutions to stochastic differential equations with stationary coefficients driven by Poisson random measures and Brownian motions.  We state an annealed convergence theorem, in which the limit  exhibits a diffusive or superdiffusive behavior, depending on the integrability properties of the Poisson random measure.\\
\noindent{\bf Keywords} : It\^o-L\'evy processes; random medium; stochastic homogenization; scaling limit; integro-differential operators; ergodicity.
\end{abstract}

\section{Introduction}\label{intro}


Consider a standard Brownian motion $\{B_t;t\geq 0\}$. It is straightforward to check that a diffusive rescaling of that process leads to the same process (in law), that is $\epsilon^{1/2}B_{t/\epsilon}$ is still a Brownian motion. This gives rise to the natural issue of determining the scaling limit of the process $X$ solution to the following Stochastic Differential Equation (SDE for short)
$$X_t=x+\int_0^tb(X_r)\,dr+\int_0^t\sigma(X_r)\,dB_r.$$
Put in other words, does the rescaled process $\epsilon^{1/2}X_{t/\epsilon}$ converge as $\epsilon\to 0$ towards a (non-standard) Brownian motion? And what does the covariations of the limiting Brownian motion look like? For several years, an extensive litterature has spread out from this topic. For a limit to exist, it is reasonable to think that the coefficients $b$ and $\sigma$ must have good averaging properties. So, the case of periodic coefficients has first been investigated (see \cite{lions,olla} for insights on the subject), and more recently, some authors have been interested in the case of stationary coefficients (see \cite{masi,kozlov,oeschlager} and many others, or \cite{sznitman} for recent issues on the topic).

On the other hand, the possible scaling limits of SDE's driven by general L{\'e}vy processes is a topic which has been poorly studied so far. This is the purpose of the following paper. More precisely, we are interested in deriving limit theorems for SDE's in stationary environments of the following form (the parameter $\omega$ stands for the randomness of the coefficients):
\begin{equation*}
X_t=x+\int_0^t ({\bfm b}+{\bfm e})(\tau_{X_{r-}}\omega)\,dr+\int_0^t \int_{\R}\gamma(\tau_{X_{r-}}\omega,z)\,\hat{N}(dr,dz)+\int_0^t{\bfm \sigma}(\tau_{X_{r-}}\omega)\,dB_r,
\end{equation*}
where $\hat{N}$ is a compensated Poisson measure. We will see that under appropriate conditions on the coefficients ${\bfm b}, {\bfm e}, \gamma, {\bfm \sigma}$, the generator of the above SDE takes the following form for sufficiently smooth functions $f$ (in a fixed environment $\omega$):
\begin{equation*}
\mathcal{L}^{\omega}f(x)=\frac{1}{2}{\bfm a}(\tau _{x}\omega)f''(x)+{\bfm b}(\tau _{x}\omega)f'(x)+\lim_{\epsilon \to 0}\int_{|z|>\epsilon}(f(x+z)-f(x)){\bfm c}(\tau _{x}\omega,z)e^{2{\bfm V}(\tau _{x}\omega)}\chi(dz),
\end{equation*}
where ${\bfm a},{\bfm b},{\bfm c},{\bfm V}$ are bounded functions of the environment.
To our knowledge, only the following papers have been devoted to deriving scaling limits of SDE's with possibly long jumps: \cite{Fra} or \cite{schwab}. Both authors consider  $ \alpha$-stable jump processes driven by periodic coefficients and treat the problem with probabilistic tools  \cite{Fra} or analytic tools \cite{schwab}. In contrast, there is an abundant litterature devoted to establishing quenched and annealed central limit theorems for SDE's driven by Poisson measures with bounded jumps. In particular, much effort has been made to derive under minimal assumptions quenched CLT's for random walks among random conductances: \cite{masi}, \cite{sido}, \cite{mathieu}.

The rest of the paper is organized as follows: in section 2, we set the notations and state the main theorem \ref{mainresult} (an annealed functional limit theorem). Section 3 and 4 are devoted to showing that the measure $\pi$ (equivalent to the original measure on the environment: see section 3 for the definition of $\pi$) is invariant for the environment seen from the particle. In section 5 and 6 are gathered some material we will need in proving the homogenization theorem (Ergodic issues and study of the correctors). In section 7 are gathered the tension estimates which are necessary to derive functional theorems in the Skorohod topology. In section 8, we give the proof of the main theorem \ref{mainresult}. Finally, in the appendices are gathered technical lemmas that are used in different places of the paper.


\section{Statements of the problem}\label{statement}

\subsection*{Random medium}
We first introduce the notion of random medium (see e.g.
\cite{jikov}) and the necessary background about random media
\begin{definition}\label{medium}
Let $(\Omega ,{\cal G},\mu )$ be a probability space and
$\left\{\tau_{x};x\in \reel\right\}$ a group of measure
preserving transformations acting ergodically on $\Omega $:

1) $\forall A\in {\cal G},\forall x\in \reel$, $\mu (\tau
_{x}A)=\mu (A)$,

2) If for any $x\in \reel$, $\tau _{x}A=A$ then $\mu (A)=0$ or
$1$,

3) For any measurable function ${\bfm g}$ on $(\Omega ,{\cal
G},\mu)$, the function $(x,\omega )\mapsto {\bfm g}(\tau_x \omega)$
is measurable on $(\reel\times\Omega ,{\cal B}(\reel)\otimes
{\cal G})$.
\end{definition}
The expectation with respect to the random medium is denoted by
${\mathbb M}$. The space of square integrable (resp. integrable, resp. essentially
bounded) functions on $(\Omega ,{\cal G},\mu )$ is denoted by
$L^2(\Omega)$ (resp. $L^1(\Omega)$, resp. $L^\infty(\Omega)$), the
usual norm by $|\, \cdot \,|_2$ (resp. $ |\, \cdot \,|_1$, resp.
$|\, \cdot \,|_\infty $) and the corresponding inner product by $(\,
\cdot\,,\, \cdot\,)_2$. The operators on $L^2(\Omega )$ defined by
$T_{x}{\bfm g}(\omega )={\bfm g}(\tau_{x}\omega )$ form a strongly
continuous group of unitary maps in $L^2(\Omega )$. Each function
${\bfm g}$ in $ L^2(\Omega )$ defines in this way a stationary
ergodic random field on $\reel $. The group possesses a
generator $D$, defined by 
\begin{equation} D{\bfm g}=\lim_{\R\ni h\to
0}h^{-1}(T_{h}{\bfm g}-{\bfm g}) \text{ if the limit exists in
the }L^2(\Omega)\text{-sense},
\end{equation}
which is closed and densely defined. We distinguish the differential operator in random medium $D$ from the usual derivative $\partial_x f $ of a function $f$ defined on $\R$.

\noindent {\bf Notations.} Recursively, we define the operators ($k\geq 1$) $D^k=D(D^{k-1})$ with domain $H^k(\Omega)=\{{\bfm f}\in H^{k-1}(\Omega); D^{k-1}{\bfm f}\in {\rm Dom}({\bfm D})=H^1(\Omega)\}$. We also define $H^\infty(\Omega)=\bigcap_{k=1}^\infty H^k(\Omega)$. 

We denote with ${\cal C}$ the dense subspace of $L^2(\Omega )$ defined by
\begin{equation*}
{\cal C}={\rm Span}\left\{{\bfm g} \star \varphi ;{\bfm g}\in
L^\infty(\Omega ),\varphi \in C^\infty _c(\reel  )\right\}
\end{equation*}
with ${\bfm g} \star \varphi(\omega )=\int_{\reel }{\bfm
g}(\tau_{x}\omega )\varphi (x)\,dx$. We point out that ${\cal
C}\subset {\rm Dom}(D)$, and $D({\bfm
g} \star \varphi)=-{\bfm g} \star
\partial  \varphi/\partial x$. This last quantity is also equal to $D{\bfm g}
\star \varphi $ if ${\bfm g}\in {\rm Dom}(D)$. $C(\Omega)$ is defined as the closure of $\mathcal{C} $ in $L^\infty(\Omega)$ with respect to the norm $|\cdot|_\infty$, whereas $C^\infty(\Omega)$ stands for the subspace of $H^\infty(\Omega)$, whose elements satisfy: ${\bfm f}\in C^\infty(\Omega)\Leftrightarrow \forall k \geq  0,\,|D^k{\bfm f} |_\infty<+\infty$. We point out that, whenever a function ${\bfm f}\in H^\infty(\Omega)$, $ \mu$
a.s. the mapping $f_\omega:x\in\R\mapsto {\bfm f}(\tau_x\omega)$  is infinitely differentiable and $\partial_x f_\omega(x)=D{\bfm f}(\tau_x\omega)$.
\subsection*{Structure of the coefficients}
We consider a so-called L\'evy measure $ \nu$, that is a  $\sigma$-finite measure  $\nu$ on $\R$ such that 
\begin{equation}\label{nu}
\int_{\R}\min(1,z^2)\nu(dz)<+\infty,\quad \nu(\{0\})=0.
\end{equation}

We introduce the coefficients $ {\bfm V},{\bfm \sigma}\in L^\infty(\Omega)$ and $ {\bfm \gamma}:\Omega\times\R\to \R$ such that  

\begin{hypothesis}\label{ellipticity}{\bf Ellipticity} We set ${\bfm a}={\bfm \sigma}^2$. There is a constant $M_{\ref{ellipticity}}>0$ such that 
$$M_{\ref{ellipticity}}^{-1}\leq {\bfm a}\leq M_{\ref{ellipticity}}.\qed$$
\end{hypothesis}
For each fixed $\omega\in \Omega$, by defining  the mapping $\gamma_\omega:z\mapsto {\bfm \gamma}(\omega,z)$, we can consider the measure $\nu\circ  \gamma_\omega^{-1}:A\subset \R\mapsto \nu(\gamma_\omega^{-1} (A))=\nu\big(\{z\in\R; {\bfm \gamma}(\omega,z)\in A\}\big)$. 

\begin{hypothesis}\label{symmetry}{\bf Symmetry of the kernel.}
We assume that the measure $\nu\circ  \gamma_\omega^{-1}$ can be rewritten as 
$$\nu\circ  \gamma_\omega^{-1}(dz)=e^{2{\bfm V}(\omega)}{\bfm c}(\omega,z)\chi(dz)$$ for some L\'evy measure $\chi$, which is symmetric (i.e. $\chi(dz)=\chi(-dz)$), and some measurable nonnegative bounded symmetric kernel ${\bfm c}$ defined on $\Omega\times \R$. The symmetry of ${\bfm c}$ means
$$\mu \text{ a.s.},\,\chi(dz) \,\text{ a.s.}, \quad {\bfm c}(\tau_z\omega,-z)={\bfm c}(\omega,z).\qed$$
 \end{hypothesis}

\begin{hypothesis}\label{regul}{\bf Regularity.} We assume the coefficients satisfy the following assumptions:

1) The coefficients ${\bfm V}$, ${\bfm \sigma}$ belong to $C^\infty(\Omega)$. In particular, we can define $${\bfm b}=\frac{1}{2}D{\bfm a}-{\bfm a}D{\bfm V} =\frac{e^{2{\bfm V} }}{2}D\big(e^{-2{\bfm V}}{\bfm a} \big)\in C^\infty(\Omega),$$ 

2) For  $\chi(dz) $-almost every $z\in\R$, the mapping $\omega\mapsto {\bfm c}(\omega,z)$ belongs to $C^\infty(\Omega)$ and, for each fixed $k\geq 1$, there exists a constant $C_k$ such that $|D^k{\bfm c}(\cdot,z)|_\infty\leq C_k$, $\chi(dz)$ a.s.

3) $\mu$ a.s., for $\nu$ almost every $|z|>1$, the mapping $x\in\R\mapsto {\bfm \gamma}(\tau_x\omega,z)$ is continuous and $\mu$ a.s., we can find a constant $C>0$ such that $\forall x,y\in\R $, $$\int_{|z|\leq 1}|{\bfm \gamma}(\tau_y\omega,z)-{\bfm \gamma}(\tau_x\omega,z)|^2\nu(dz)\leq C |y-x|^2,\,\,\,     \int_{|z|\leq 1}|{\bfm \gamma}(\tau_x\omega,z)|^2\nu(dz)\leq C(1+|x|^2).$$

4) The limit $${\bfm e}(\omega)=\lim_{\alpha\to 0}\int_{\alpha\leq |{\bfm \gamma}(\omega,z)|}{\bfm \gamma}(\omega,z)\one_{|z|\leq1}\nu(dz)$$ exists in the $ L^2(\Omega)$ sense and defines bounded Lipschitzian function, that is (for some constant $M_{\ref{regul}}\geq 0$), $|{\bfm e}|_\infty\leq M_{\ref{regul}}$ and $ \mu$  a.s., $\forall x,y\in\R$, $ |{\bfm e}(\tau_y\omega)-{\bfm e}(\tau_x\omega)|\leq M_{\ref{regul}}|x-y|$. Furthermore, there is a positive constant $S$ such that $\sup_{|z|\leq 1}|{\bfm \gamma}(\cdot,z)|_{\infty}\leq S$
\qed
 \end{hypothesis}

\begin{hypothesis}\label{fouriersc}{\bf Convergence rate.} 
We assume either of the following conditions holds:
\begin{enumerate}
\item {\bf (pure jump scaling)} In the case $\int_{\R}z^2\chi(dz)=+\infty$, we assume that there are a function $\delta:]0;+\infty[\to]0;+\infty[$ satisfying $\lim_{\epsilon\to 0}\delta(\epsilon)=0 $, a non-zero random function $ {\bfm \theta}:\{-1;1\}\to L^\infty(\Omega)$ and a L\'evy measure $\mathcal{H}$ on $\R$ such that
\begin{equation}\label{cr}
\lim_{\epsilon\to 0}\M\big[|\epsilon^{-1}\int_{\R}g(\delta(\epsilon)z) {\bfm c}(\cdot,z)\chi(dz)-\int_{\R}{\bfm \theta}(\cdot,{\rm sign}(z))g(z)\mathcal{H}(dz)|\big]=0
\end{equation}
for each function $g=\one_{[a,b]}$, with $a<b$ and $0\not\in [a,b]$. Throughout the paper, the random measure ${\bfm \theta}(\omega,{\rm sign}(z))\mathcal{H}(dz)$ will be called the limit measure. 

We further require the quantity  $\epsilon^{-1}\delta(\epsilon)^2\int_{\delta(\epsilon)|z|\leq \alpha}z^2\chi(dz)$ to be converging towards $0$ as $\alpha\downarrow 0$, uniformly with respect to $\epsilon$.

We point out that, necessarily in that case, $\lim_{\epsilon\to 0}\delta(\epsilon)^2/\epsilon=0$. 
 
\item {\bf (diffusive scaling)} In the case $\int_{\R}z^2\chi(dz)<+\infty$, we set $\delta(\epsilon)=\epsilon^{1/2}$. 
\end{enumerate}
\end{hypothesis}

\begin{remark}
Let us make a few comments about our assumptions. \ref{regul}.3 and \ref{regul}.4 are only technical assumptions to ensure existence and uniqueness of a solution to SDE \eqref{SDE} below, whereas \ref{regul}.1, \ref{regul}.2 make the resolvent operator associated to \eqref{SDE} regularizing enough. Assumptions \ref{symmetry} and \ref{fouriersc} are closely related to the scaling properties of \eqref{SDE}. In particular, \ref{fouriersc} states that the jump measure possesses good averaging properties.
\end{remark}
Even if it means adding to ${\bfm V}$ a renormalization constant (this does not change the drift ${\bfm b}$ and the jump coefficients ${\bfm \gamma} $ and $\nu$), we consider the probability measure $d\pi=e^{-2{\bfm V}}\,d\mu$ on $(\Omega,\mathcal{G})$, and we denote by $\M_\pi$ the expectation w.r.t. this probability measure.

\subsection*{Jump-diffusion processes in random medium}
 We suppose that we are
given a complete probability space  $(\Omega',{\cal F},\P)$ with a
right-continuous increasing family of complete sub $\sigma$-fields $({\cal F}_t)_t$
of ${\cal F}$, a ${\cal F}_t$-adapted Brownian motion $ \{B_t;t\geq 0\}$ and ${\cal F}_t$-adapted Poisson random measure $N(dt,dz)$ with intensity $\nu$. $\tilde{N}(dt,dz)=N(dt,dz)-\nu(dz)dt$ denotes the corresponding compensated random measure and $\hat{N}(dt,dz)$ the truncated compensated random measure $N(dt,dz)-\one_{|z|\leq 1}\nu(dz)dt$. We further assume that the Brownian motion, the L\'evy process and the random medium are independent. 

For each fixed $ \omega\in\Omega$, Assumptions \ref{regul}.3 and \ref{regul}.4 are enough to ensure  existence and pathwise uniqueness of a ${\cal F}_t$-adapted process $X$ (see \cite[Ch.6, Sect.2]{applebaum}) solution to the following SDE
\begin{equation}\label{SDE}
\begin{split}
X_t=x+\int_0^t ({\bfm b}+{\bfm e})(\tau_{X_{r-}}\omega)\,dr+\int_0^t \int_{\R}\gamma(\tau_{X_{r-}}\omega,z)\,d\hat{N}(dr,dz)+\int_0^t{\bfm \sigma}(\tau_{X_{r-}}\omega)\,dB_r.
\end{split}
\end{equation}

\subsection*{Main result}
We denote with $C(\R_+;\R)$ the space of continuous $\R$-valued functions on $[0;+\infty[$, endowed with the topology of uniform convergence on compact intervals and with $D(\R_+;\R)$  the space of right-continuous $\R$-valued functions with left limits, endowed with the Skorohod topology, cf \cite{ethier}. We claim
\begin{theorem}\label{mainresult}
{\bf 1) Pure jump scaling:} In the case $\int_\R z^2\chi(z)dz=+\infty$, in $\mu$ probability, the rescaled process $\delta(\epsilon)X_{\cdot/\epsilon}$, starting from $0\in \R$, converges in law towards a L\'evy process with L\'evy symbol
$$\int_{\R}(e^{iuz}-1-i u z\one_{|z|\leq 1})\M[{\bfm \theta}(\cdot,{\rm sign}(z) )]\mathcal{H}dz) $$ in the Skorohod topology.

{\bf 2) Diffusive scaling:} In the case $\int_\R z^2\chi(z)dz<+\infty$, in $\mu$ probability, the process $X$, starting from $x\in \R$, converges in law in the Skorohod topology towards a non standard centered Brownian motion with variance  $A$ given by (see Section \ref{sec:correctors} for the definition of ${\bfm \xi}$ and $ {\bfm \zeta}$)
\begin{equation}\label{form:A}
A= \M\big[{\bfm a}(1+{\bfm \xi})^2e^{-2{\bfm V}}+\int_\R(z+{\bfm \zeta}(\cdot,z))^2{\bfm c}(\cdot,z)\chi(dz)\big]
\end{equation}
\end{theorem}

\begin{remark}
Actually, by adapting the proof of \cite[Section 2.7]{olla}, we can prove that $A$ is given by the variational formula
\begin{equation}\label{variational}
A=\inf_{{\bfm \varphi}\in \mathcal{C}} \M\big[{\bfm a}(1+D{\bfm \varphi})^2e^{-2{\bfm V}}+\int_\R(z+T_z{\bfm \varphi}-{\bfm \varphi})^2{\bfm c}(\cdot,z)\chi(dz)\big],
\end{equation}
from which lower and upper bounds for $A$ can be obtained. In particular, $A$ is nondegenerate (because ${\bfm a}$ is).
\end{remark}
\begin{remark}
We stress that our result is stated in dimension 1 but our proofs straightforwardly extend to higher dimensions, though it might be notationally more challenging.
\end{remark}

\subsection*{Applications}
 Suppose the jump rate ${\bfm c}(\omega,z)\chi(dz)$ is known  It gives rise to the issue of determining a coefficient ${\bfm \gamma}$ and a measure $\nu$ satisfying Assumptions \ref{symmetry} and \ref{regul}. In most classical situations, the  following lemma is helpful to construct such a  ${\bfm \gamma}$:
 
\begin{lemma}{\bf Generic construction of a coefficient ${\bfm \gamma}$ and measure $\nu$ associated to a prescribed jump rate of the form $\mathbf{{\bfm c}(\omega,z)\chi(dz)}$:}\label{constructiongamma} Suppose we are given ${\bfm c}:\Omega\times \R\to ]0,+\infty[$ and a  strictly positive even function $\chi:\R\to ]0,+\infty[$, bounded on the compact subsets of $\R\setminus \{0\}$, satisfying:

1) $\chi(z)dz$ is a L\'evy measure such that (for some positive constant $M'$)
$$\int_{0}^{+\infty}\chi(z)\,dz=+\infty,\quad \text{and} \quad \forall z\in]0,1],\quad \int_z^{+\infty}\chi(r)\,dr\leq M'\chi(z)z ,$$

2) for some  constants $0<m\leq M$, we have $m\leq {\bfm c}(\omega,z)\leq M $. 

3) $\forall z \in\R$, $ {\bfm c}(\cdot,z)\in C^\infty(\Omega)$ and $\forall k\geq 1$, $\exists C_k\geq 0$, $|D^k{\bfm c}(\cdot,z)|_\infty\leq C_k$.

Under the above assumptions, by setting $${\bfm c}^s(\omega,z)=\frac{1}{2}\big({\bfm c}(\tau_z\omega,-z)+{\bfm c}(\omega,z)\big),$$  we define a symmetric kernel fitting all the conditions required in Assumption \ref{symmetry} and \ref{regul}. Moreover, we can find a coefficient ${\bfm \gamma}:\Omega\times \R\to \R$ and a L\'evy measure $\nu$ fitting the regularity conditions of points 3) and 4) of Assumption \ref{regul}, satisying $|{\bfm \gamma}(\omega,z)|\leq|z|$ and such that the measures $\nu\circ {\bfm \gamma}^{-1}(\omega,\cdot)$ and ${\bfm c}^s(\omega,z)\chi(dz)$ coincide. 
\end{lemma}
\begin{remark}
For instance, for any $\alpha\in ]0,2[$ and $\beta\in\R$, the L\'evy measures $$\chi(z)=|z|^{-1-\alpha},|z|^{-1-\alpha}(\ln(1+ |z|))^\beta,e^{-|z|} |z|^{-1-\alpha},\cdots$$ (and many others) suit. 
\end{remark}
\begin{remark} With a little care, we can adapt the proof of Lemma \ref{constructiongamma} to treat the cases  $\int_0^{+\infty}\chi(z)\,dz<+\infty $ or $\chi(z)\not=\chi(-z)$. What really matters in the proof is the condition $\chi(z)>0$.  
\end{remark}
\begin{lemma}{\bf Case of pure jump scaling.}\label{lemmapjs} In the  case $\int_{\R}z^2\chi(z)dz=+\infty$, suppose the following conditions hold:
 
4)  for some non-zero functions ${\bfm \theta}: \{-1;1\}\to L^\infty(\Omega)$ $$\lim_{z\to \pm\infty}\M[|{\bfm c}(\omega,z)-{\bfm \theta}(\omega,\pm 1)|]=0.$$ 

5) there is a L\'evy measure $\mathcal{H}(dz)$ such that, for any function $g=\one_{[a,b]}$ (with $0\not\in [a,b]$), we have
\begin{equation*}\label{condfour}
\forall u\in\R,\quad \lim_{\epsilon\to 0}\frac{1}{\epsilon}\int_{\R}g(\delta(\epsilon)z)e^{izu}\chi(dz)=\one_{u=0}\int_{\R}g(z)\mathcal{H}(dz).
\end{equation*}

6) The quantity  $\epsilon^{-1}\delta(\epsilon)^2\int_{\delta(\epsilon)|z|\leq \alpha}z^2\chi(dz)$ converges towards $0$ as $\alpha\downarrow 0$, uniformly with respect to $\epsilon$. 

Under the avove assumptions, Assumption \ref{fouriersc} is satisfied with convergence rate $\delta(\epsilon)$ and limit measure  $\Big({\bfm \theta}(\omega,{\rm sign}(z))+\M\big[{\bfm \theta}(\omega,-{\rm sign}(z))\big]\Big)\mathcal{H}(dz) $.
\end{lemma}

\begin{remark}
When the measure $\chi(dz)$ is of the type $\chi(dz)=\frac{1}{|z|^{1+\alpha}}dz$ for some $\alpha\in]0,2[$, point 5) is particularly easy to check since it results from the Riemann-Lebesgue theorem after choosing $\delta(\epsilon)=\epsilon^{1/\alpha}$ and making the change of variables $ \epsilon^{1/\alpha}z=y$.
\end{remark}

Concerning Examples \ref{alphast}-\ref{attract} below, the superscript $s$ of a symmetric kernel ${\bfm c}^s$ means that ${\bfm c}^s$ is constructed as prescribed by the generic construction (Lemma \ref{constructiongamma}) from a reference function ${\bfm c}$.  Moreover, all the considered L\'evy measures $\chi $ satisfy the conditions of  the generic construction and of Lemma \ref{lemmapjs}. The reference function ${\bfm c}$ is assumed to be converging towards a function $ {\bfm \theta}:\{-1;1\}\to L^\infty(\Omega)$. It is thus convenient to define
$$\Theta(\omega,z)={\bfm \theta}(\omega,{\rm sign}(z))+\M[{\bfm \theta}(\omega,-{\rm sign}(z))].$$ 
To sum up, in examples  \ref{alphast}-\ref{attract} below, given a triple $({\bfm \sigma}, {\bfm c}^s,\chi)$, we can construct the corresponding coefficients ${\bfm \gamma}$ and $\nu$, define the process $X$ solution of \eqref{SDE} and apply Theorem \ref{mainresult}. So we won't specify these points anymore. We just precise, in each case, what the limit measure and convergence rate look like. 

\begin{example}{{\bf $\mathbf{\alpha}$-stable kernels.}}\label{alphast}
We consider the kernel  $$\frac{{\bfm c}^s(\omega,z)}{|z|^{1+\alpha}}dz,\quad  0<\alpha<2.$$  
 The convergence rate is given by $\delta(\epsilon)=\epsilon^{1/\alpha}$ and the limit measure by $\frac{\Theta(\omega,z)}{|z|^{1+\alpha}}dz.$
\end{example}

\begin{example}{{\bf Multi-stable kernels.}} 
Given a parametrized family  $({\bfm c}^s(\cdot,z,\alpha))_{\alpha_1\leq \alpha\leq \alpha_2}$ ($\alpha_1,\alpha_2\in]0,2[ $), we  are interested in the kernel
$$\int_{\alpha_1}^{\alpha_2}\frac{{\bfm c}^s(\omega,z,\alpha)}{|z|^{1+\alpha}}\,d\alpha.$$ The coefficient ${\bfm \gamma}$ can be constructed from the L\'evy measure $\chi(dz)=\int_{\alpha_1}^{\alpha_2}\frac{1}{|z|^{1+\alpha}}\,d\alpha$. Then Assumption \ref{fouriersc} holds with convergence rate $\delta(\epsilon)$ given by the implicit relation $\lim_{\epsilon\to 0}\delta(\epsilon)^{\alpha_1}(-\epsilon\ln(\delta(\epsilon)))^{-1}=1$. The limit measure matches $\frac{\Theta(\omega,z,\alpha_1)}{|z|^{1+\alpha_1}}dz.$
\end{example}

\begin{example}{{\bf kernels attracted by stable kernels.}}\label{attract}
We can generalize Example \ref{alphast} as follows. Given $0<\alpha<2$ and a bounded function $l:]0;+\infty[\to \R$ such that $\lim_{r\to +\infty}l(r)=0$, we define  $h(z)=\exp(\int_0^{|z|}l(r)r^{-1}\,dr)$ and $$\chi(z)=\frac{h(z)}{|z|^{1+\alpha}}.$$
Without giving further details, the reader may check criterion \eqref{condfour} with $\mathcal{H}(dz)=\frac{1}{|z|^{1+\alpha}}$ and convergence rate implicitly given by the (asymptotic) relation: $\delta(\epsilon)^{\alpha}h(1/\delta(\epsilon))/\epsilon\to 1$ as $\epsilon\to 0$.

The most famous examples are given by ($\beta_1,\beta_2,\dots\in\R$) $h(z)=\big(\ln(|z|+1)\big)^{\beta_1}, \big(\ln(|z|+1)\big)^{\beta_1}\big(\ln\big(1+\ln(|z|+1)\big)^{\beta_2},\dots $ and so on. Now, looking at the kernel ${\bfm c}^s(\omega,z)\chi(z)dz$, the limit measure is given by $\frac{{\bfm \Theta}(\omega,z)}{|z|^{1+\alpha}}dz.$
\end{example}

Now, we investigate the situation when the kernel  ${\bfm c}(\omega,z)\chi(dz) $ corresponds to that of a random walk.
\begin{example}{{\bf Symmetric random walk among random conductances.}}\label{rwrc}
Consider a smooth random variable ${\bfm W}:\Omega\to [a,b]$ for some constants $b>a>0$. Define $\chi=\delta_1+\delta_{-1}$ ($\delta_a$ denotes the Dirac mass at point $a\in\R$), ${\bfm c}(\cdot,1)=T_{1/2}{\bfm W}$ and ${\bfm c}(\cdot,-1)=T_{-1/2}{\bfm W}$, and ${\bfm V}=-(1/2)\ln(T_{-1/2}{\bfm W}+T_{1/2}{\bfm W})$. The kernel
$$e^{2{\bfm V}}{\bfm c}(\cdot,z)\chi(dz)=\frac{T_{-1/2}{\bfm W}}{T_{-1/2}{\bfm W}+T_{1/2}{\bfm W}}\delta_{-1}+\frac{T_{1/2}{\bfm W}}{T_{-1/2}{\bfm W}+T_{1/2}{\bfm W}}\delta_{1}$$ corresponds to that of a random walk among random conductances. Set $$\nu(dz)=\one_{[0,1]}(z)dz,\quad {\bfm \gamma}(\omega,z)=\left\{\begin{array}{ll}1& \text{ if }z\leq {\bfm c}(\omega,1)e^{2{\bfm V}},\\
-1& \text{ if }z>{\bfm c}(\omega,1)e^{2{\bfm V}}.\end{array}\right.$$
Clearly, the measures $\nu\circ \gamma_\omega^{-1}$ and $e^{2{\bfm V}}{\bfm c}(\cdot,z)\chi(dz)$ coincide. The reader may easily check that the regularity conditions of Assumption \ref{regul} are satisfied. Moreover, we are clearly in the situation of diffusive scaling and, in case ${\bfm a}=1$, Theorem \ref{mainresult} ensures that a  mixed Brownian motion/random walk among random conductances behaves like a Brownian motion with effective diffusivity
$$A=\inf_{{\bfm \varphi}\in\mathcal{C}}\M_\pi\big[(1+D{\bfm \varphi})^2e^{-2{\bfm V}}+{\bfm c}(\cdot,1)(1+T_1{\bfm \varphi}-{\bfm \varphi})^2+{\bfm c}(\cdot,-1)(-1+T_{-1}{\bfm \varphi}-{\bfm \varphi})^2\big]. $$
\end{example}

\section{Dirichlet forms in random medium}\label{sec:dirichlet}

For the sake of readibility, the proofs of this section are gathered in Appendix \ref{app:dirichlet} and may be omitted upon the first reading.

We can then equip the space $L^2(\Omega)$ with the inner product $({\bfm \varphi},{\bfm \psi})_\pi=\M[{\bfm \varphi}{\bfm \psi}e^{-2{\bfm V}}]$, and denote by $|\cdot|_\pi$ the associated norm. Since ${\bfm V}$ is bounded, both inner products $(\cdot,\cdot)_2$ and $ (\cdot,\cdot)_\pi$ are equivalent on $L^2(\Omega)$.


Let us define on ${\cal C}\times {\cal C} $ the following bilinear forms (with $\lambda>0$) 
\begin{equation}\label{dirichletf}
\begin{split}
B^d({\bfm \varphi},{\bfm \psi})&=\frac{1}{2}({\bfm a}D{\bfm \varphi},D{\bfm \psi})_\pi,\quad B^{j}({\bfm \varphi},{\bfm \psi})=\frac{1}{2}\M\int_{\R}(T_z{\bfm \varphi}-{\bfm \varphi})(T_z{\bfm \psi}-{\bfm \psi}){\bfm c}(\cdot,z)\chi(dz),\\
B^s({\bfm \varphi},{\bfm \psi})&=B^d({\bfm \varphi},{\bfm \psi})+B^{j}({\bfm \varphi},{\bfm \psi}),\quad
B^s_\lambda({\bfm \varphi},{\bfm \psi})=\lambda({\bfm \varphi},{\bfm \psi})_\pi+B^s({\bfm \varphi},{\bfm \psi}).
\end{split}
\end{equation}
We can thus consider on ${\cal C}\times {\cal C}  $ the inner product $B^s_\lambda $ and the closure $\H$ of ${\cal C} $ w.r.t. the associated norm (note that the definition of $\H$ does not depend on $\lambda>0$ since the corresponding norms are equivalent). From now on, our purpose is to construct a self-adjoint operator associated to $B^s$ and to derive its regularizing properties.
 
The following construction follows \cite[Ch. 3, Sect. 3]{fuku} (or \cite[Ch. 1, Sect. 2]{ma}), to which the reader is referred for further details. For any $\lambda>0$, $B_\lambda$ is clearly continuous on ${\cal C}\times {\cal C}  $ so that it continuously extends to $\H\times \H$ (the extension is still denoted $B_\lambda$). Moreover, $B_\lambda^s$ is coercive. It thus defines a resolvent operator $G_\lambda: L^2(\Omega)\to \H$, which is one-to-one and continuous. We define the unbounded operator ${\bfm L}=\lambda-G_\lambda^{-1}$ on $L^2(\Omega)$ with domain ${\rm Dom}({\bfm L})=G_\lambda(L^2(\Omega))$. This definition does not depend on $\lambda>0$. More precisely, a function
${\bfm \varphi}\in \H$ belongs to ${\rm Dom}({\bfm L})$ if and only if the map
${\bfm \psi}\in \H\mapsto B^s_\lambda({\bfm \varphi},{\bfm \psi})$ is $L^2(\Omega)$
continuous. In this case, we can find ${\bfm f}\in L^2(\Omega)$ such
that $ B_\lambda({\bfm \varphi},\cdot)=({\bfm f},\cdot)_\pi$. Then ${\bfm L}{\bfm \varphi}$
exactly matches $\lambda{\bfm \varphi}-{\bfm f}$. We point out that the unbounded operator
${\bfm L}$ is closed, densely defined and seld-adjoint. We further stress that the weak form of the resolvent equation $\lambda G_\lambda{\bfm f}-{\bfm L}G_\lambda{\bfm f}={\bfm f} $ reads: $\forall {\bfm \psi}\in\H$
\begin{align}\label{eqbase}
\lambda(G_\lambda{\bfm f},{\bfm \psi})_\pi+&\frac{1}{2}({\bfm a}DG_\lambda{\bfm f},D{\bfm \psi})_\pi+\frac{1}{2}\M\int_{\R}(T_zG_\lambda{\bfm f}-G_\lambda{\bfm f})(T_z{\bfm \psi}-{\bfm \psi}){\bfm c}(\cdot,z)\chi(dz)\\&= ({\bfm f},{\bfm \psi})_\pi.\nonumber
\end{align}

For sufficiently smooth functions, ${\bfm L}$ can be easily identified:
\begin{lemma}\label{genito}
 Let ${\bfm \varphi}\in H^\infty(\Omega)$. Then ${\bfm \varphi}\in{\rm Dom}({\bfm L})$
\begin{align}\label{generator}
 {\bfm L}{\bfm \varphi}&=\frac{1}{2}{\bfm a}D^2{\bfm \varphi}+({\bfm b}+{\bfm e})D{\bfm \varphi}+\int_{\R}\big({\bfm \varphi}(\tau_{{\bfm \gamma}(\omega,z)}\omega)- {\bfm \varphi}(\omega)-{\bfm \gamma}(\omega,z)\one_{\{|z|\leq 1\}}D{\bfm \varphi}(\omega)\big)\,\nu(dz).
\end{align}
\end{lemma}

We now investigate the regularizing properties of the resolvent operator $G_\lambda$.
\begin{proposition}\label{core}
For each $\lambda>0$, the resolvent operator $G_\lambda$ maps $L^2$ into $H^2(\Omega)$, and $H^m(\Omega)$ into $H^{m+2}(\Omega)$ for any $m\geq 1$. In particular ${\rm Dom}({\bfm L}^m)=H^{2m}(\Omega) $.
\end{proposition}

The operator ${\bfm L}$ is self-adjoint. Thus it generates a strongly continuous contraction semi-group $(P_t)_t$ of self-adjoint operators. Each operator $P_t$ ($t>0$) maps $L^2(\Omega)$ into ${\rm Dom}({\bfm L})=G_\lambda(L^2(\Omega))\subset H^2(\Omega)$. More precisely, combining Hille-Yosida's theorem with Proposition \ref{core}, we get the following estimates:
 \begin{align}\label{reg:semigroup}
 {\bfm f}\in L^2(\Omega)&\Rightarrow t\mapsto P_t{\bfm f} \in C([0;+\infty[;L^2(\Omega))\cap C^\infty(]0;+\infty[;H^\infty(\Omega)),\\
{\bfm f}\in H^\infty(\Omega)&\Rightarrow t\mapsto P_t{\bfm f} \in C^\infty([0;+\infty[;H^\infty(\Omega)).
\end{align}
where, given an interval $I\subset\R$, $C(I;L^2(\Omega))$ (resp. $C^\infty(I;H^\infty(\Omega))$) stands for the space of continuous functions from $I$ to $L^2(\Omega)$ (resp. infinitely differentiable functions from $I$ to $ H^\infty(\Omega)$). Moreover, we can prove
\begin{proposition}\label{submarkov}
The semi-group $(P_t)_t$ is sub-Markovian. Put in other words, for any ${\bfm f}\in L^2(\Omega)$ such that $0\leq {\bfm f}\leq 1$ $\mu$ a.s., we have $0\leq P_t{\bfm f}\leq 1$ $\mu$ a.s. for any $t>0$. In particular, $P_t:L^\infty(\Omega)\to L^\infty(\Omega)$ and $G_\lambda: L^\infty(\Omega)\to L^\infty(\Omega)$ are continuous.
\end{proposition}

\section{Environment as seen from the particle}\label{particle}

In what follows, $X$ denotes the solution of \eqref{SDE} starting from $0$. Let us consider a bounded function ${\bfm \varphi}\in C^\infty([0,+\infty[; H^\infty( \Omega))$. In particular, $\mu$ a.s., the mapping $(t,x)\mapsto {\bfm \varphi}(t,\tau_x\omega)$ belongs to $C^\infty([0,+\infty[\times\R)$ and is bounded.

We can thus apply the It\^o formula (see \cite[Ch. II, Th. 32]{protter} or \cite[Ch. III]{applebaum}): $\mu$ a.s.
\begin{align*}
{\bfm \varphi}&(t,\tau_{X_t}\omega)={\bfm \varphi}(0,\omega)+\int_0^t\big(\partial_t{\bfm \varphi}+\frac{1}{2}{\bfm a}D^2{\bfm \varphi}+{\bfm b}D{\bfm \varphi}+{\bfm e}D{\bfm \varphi}\big)(r,\tau_{X_{r-}}\omega)\,dr\\ &+\int_0^tD{\bfm \varphi}{\bfm \sigma}(r,\tau_{X_{r-}}\omega)\,dB_r+\int_0^t({\bfm \varphi}(r,\tau_{X_{r-}+{\bfm \gamma}(\tau_{X_{r-}}\omega,z)}\omega)-{\bfm \varphi}(r,\tau_{X_{r-}}\omega))\,\tilde{N}(dr,dz)\\&\!\!+\!\int_0^t\!\!\int_{\R}\!\!\big({\bfm \varphi}(r,\!\tau_{X_{r\!-}\!+{\bfm \gamma}(\!\tau_{X_{\!r\!-}}\!\omega,z)}\omega\!)\!-\!{\bfm \varphi}(r,\!\tau_{X_{r\!-}}\!\omega\!)\!-\!{\bfm \gamma}(\!\tau_{X_{r\!-}}\!\omega,z\!)\one_{\{|z|\!\leq \!1\}}\!D{\bfm \varphi}(r,\tau_{X_{r\!-}}\omega)\big)\,\nu(dz)dr.
\end{align*}
It is thus natural to investigate the properties of the $\Omega$-valued process $Y_t(\omega)=\tau_{X_{t}}\omega$, which is Markovian as a consequence of \cite[Th. 6.4.6]{applebaum} and its
generator coincides on ${\cal C}$ with ${\bfm L}$ from the above computations.

\begin{proposition}
For each function $ {\bfm f}\in C(\Omega)$, we have $P_t{\bfm f}(\omega)=\E[{\bfm f}(Y_t(\omega))]$ $\mu$ a.s. As a consequence, $\forall {\bfm f}\in C(\Omega)$, $\M_\pi[\E[{\bfm f}(Y_t(\omega))]]=\M_\pi[\E[{\bfm f}(Y_{t-}(\omega))]]= \M_\pi[{\bfm f}]$.
\end{proposition}

\noindent {\bf Proof.}  From \eqref{reg:semigroup}, given ${\bfm \varphi}\in H^\infty(\Omega)\cap L^\infty(\Omega)$ and $t>0$, the mapping $(s,\omega)\mapsto P_{t-s}{\bfm \varphi}$ belongs to $ C^\infty([0,t]; H^\infty( \Omega))$ and is bounded (cf Prop \ref{submarkov}). We can thus apply the above It\^o formula between $ 0$ and $t$, which reads (use $\partial_tP_t{\bfm \varphi}={\bfm L} P_t{\bfm \varphi}$) $\mu$ a.s.:
\begin{align}\label{tim}
{\bfm \varphi}&(\tau_{X_t}\omega)=P_t{\bfm \varphi}(\omega)+\int_0^tDP_{t-r}{\bfm \varphi}{\bfm \sigma}(\tau_{X_{r-}}\omega)\,dB_r\\&+\int_0^t(P_{t-r}{\bfm \varphi}(\tau_{X_{r-}+{\bfm \gamma}(\tau_{X_{r-}}\omega,z)}\omega)-P_{t-r}{\bfm \varphi}(\tau_{X_{r-}}\omega))\,\tilde{N}(dr,dz)\nonumber
\end{align}
We remind the reader that $\mu$ a.s., $\P($ X $\text{ is c\`ad-l\`ag on }[0,t])=1$. Hence $\P(\sup_{0\leq s \leq t}|X_s|<+\infty)=1$. We deduce that the sequence of stopping times $ S_n=\inf\{s\geq 0; |X_s|>n\}$ satisfies: $\mu$ a.s., $\P$ a.s. $S_n\to +\infty$ as $n\to \infty$. By replacing $t$ by $t\wedge S_n$ (i.e. $\min(t,S_n) $) in \eqref{tim} and by taking the expectation, the martingale terms vanish and we get
$$\E[{\bfm \varphi}(\tau_{X_{t\wedge S_n}}\omega)] =\E[P_{t\wedge S_n}{\bfm \varphi}(\omega)].$$
Using the boundedness of ${\bfm \varphi}$ and $P_t{\bfm \varphi}$, we can pass to the limit as $n\to \infty$ in the above equality  to prove $P_t {\bfm \varphi}(\omega)=\E[{\bfm \varphi}(\tau_{X_{t}}\omega)]=\E[{\bfm \varphi}(Y_t(\omega))]$. In case ${\bfm f}\in C(\Omega)$, we can find a sequence $({\bfm \varphi}_n)_n \in C^\infty(\Omega)$ converging towards ${\bfm f}$ in $L^\infty(\Omega)$-norm (for instance $({\bfm f}\star \rho_n)_n$ for some regularizing sequence $(\rho_n)_n\subset C^\infty_c(\R)$). We complete the proof by passing to the limit in the relation $P_t {\bfm \varphi}_n(\omega)=\E[{\bfm \varphi}_n(\tau_{X_{t}}\omega)]=\E[{\bfm \varphi}_n(Y_t(\omega))]$.
\qed

\begin{corollary}
The measure $\pi$ is invariant for the Markov process $Y$. 
\end{corollary}
\section{Ergodic problems}
 
This section is devoted to the study of the asymptotic properties of  the process $Y$. As illustrated below, this is deeply connected to the behaviour of the resolvent $G_\lambda$ when $\lambda$ goes to $0$.

\begin{theorem}{\bf Ergodic theorem I.}\label{ergth}
For any ${\bfm f}\in L^1(\Omega)$, the following convergence holds
$$\lim_{t\to \infty}\M_\pi\E\Big[\big|\frac{1}{t}\int_0^t{\bfm f}(\tau_{X_{r-}}\omega)\,dr-\M_\pi[{\bfm f}]\big|\Big]=0. $$
\end{theorem}

\noindent {\bf Proof.} This is nothing but the ergodic theorem for stationary Markov processes (see \cite{daprato}). However, it remains to check that the measure is ergodic for the process $Y$, that is for any function ${\bfm f}\in L^2(\Omega)$ satisfying $P_t{\bfm f}={\bfm f}$ $\pi$ a.s. for any $\forall t>0$, then ${\bfm f}$ is constant $\pi$ a.s.. Such a function ${\bfm f}$ necessarily belongs to ${\rm Dom }({\bfm L})$  and satisfies $ {\bfm L}{\bfm f}=0$. Hence $ {\bfm f}\in \H $ and $B^s({\bfm f},{\bfm f})=0$. In particular $D{\bfm f}=0$ (because of Assumption \ref{ellipticity}), i.e. ${\bfm f}$ is constant $\mu$ almost surely. Since $\mu$ and $\pi$ are equivalent, we complete the proof. \qed

\begin{remark}\label{cvz2}
Under Assumption \ref{fouriersc} (case of pure jump scaling), it is a simple exercise to show that
the convergence
$$\lim_{\epsilon\to 0}\M\big[|\epsilon^{-1}\int_{\R}g(\delta(\epsilon)z) {\bfm c}(\cdot,z)\chi(dz)-\int_{\R}{\bfm \theta}(\cdot,{\rm sign}(z))g(z)\mathcal{H}(dz)|\big]=0$$ actually holds for any piecewise continuous function $g$ satisfying: $\forall z \in\R$,  $|g(z)|\leq M\min(z^2,1)$ for some positive constant $M$. 
\end{remark}

\begin{corollary}\label{coro_ergodic}
{\bf 1) Case of pure jump scaling:} Consider a piecewise continuous function $g$ satisfying: $\forall z \in\R$,  $|g(z)|\leq M\min(z^2,1)$ for some positive constant $M$, and define
$$ {\bfm G}^c_\epsilon(\omega)=\frac{1}{\epsilon}\int_{\R}g(\delta(\epsilon)z){\bfm c}(\omega,z)e^{2{\bfm V}(\omega)}\chi(dz),\quad {\bfm G}^\mathcal{H}(\omega)= \int_{\R}g(z){\bfm \theta}(\omega,{\rm sign}(z))e^{2{\bfm V}(\omega)}\mathcal{H}(dz).$$
The following convergence holds 
$$\lim_{\epsilon\to 0}\M_\pi\E\Big[\big|\epsilon\int_0^{\frac{t}{\epsilon}} {\bfm G}^c_\epsilon(\tau_{X_{r-}}\omega)dr- t\M_\pi[{\bfm G}^\mathcal{H}]\big|\Big]=0.$$ 
{\bf 2) Case of diffusive scaling:} Consider a measurable function ${\bfm g}:\Omega\times\R\to\R$ such that $\M\int_\R|{\bfm g}(\cdot,z)|{\bfm c}(\cdot,z)\chi(dz)<+\infty$, and define $${\bfm G}(\omega)=\int_\R{\bfm g}(\omega,z){\bfm c}(\omega,z)\chi(dz).$$ Then we have
$$\lim_{\epsilon\to 0}\M_\pi\E\Big[\big|\epsilon\int_0^{t/\epsilon}{\bfm G}(\tau_{X_{r-}}\omega)dr- t\M_\pi[{\bfm G}]\big|\Big]=0.$$
\end{corollary}

\noindent {\bf Proof.}  1) Case of pure jump scaling: Define 
$$ {\bfm G}^c_\epsilon(\omega)=\frac{1}{\epsilon}\int_{\R}g(\delta(\epsilon)z){\bfm c}(\omega,z)e^{2{\bfm V}(\omega)}\chi(dz),\quad {\bfm G}^\mathcal{H}(\omega)= \int_{\R}g(z){\bfm \theta}(\omega,{\rm sign}(z))e^{2{\bfm V}(\omega)}\mathcal{H}(dz).$$
By using the invariance of the measure $\pi$ for the process $Y$, we have:
\begin{align*}
\M_\pi\Big[\E|&\epsilon\int_0^{t/\epsilon}{\bfm G}^c_\epsilon(Y_{r-}(\omega))dr-\epsilon\int_0^{t/\epsilon}{\bfm G}^\mathcal{H}(Y_{r-}(\omega))dr|\Big]\\&\leq \epsilon\int_0^{t/\epsilon}\M_\pi\Big[\E|{\bfm G}^c_\epsilon(Y_{r-}(\omega))-{\bfm G}^\mathcal{H}(Y_{r-}(\omega))|\Big]dr\\
&= t\M_\pi\big[\big|{\bfm G}^c_\epsilon-{\bfm G}^\mathcal{H}\big|\big],
\end{align*}
and this latter quantity tends to $0$ as $\epsilon\to 0$ in virtue of Remark \ref{cvz2}. Since  $ {\bfm G}^\mathcal{H}$ belongs to $L^1(\Omega)$, Theorem \ref{ergth} establishes that $ \epsilon\int_0^{t/\epsilon}{\bfm G}^\mathcal{H}(Y_{r-}(\omega))dr$ converges to $ t\M_\pi[{\bfm G}^\mathcal{H}]$ as $\epsilon\to 0$. We complete the proof in that case.

2) Case of diffusive scaling: since ${\bfm G}\in L^1(\Omega)$, this is a direct consequence of Theorem \ref{ergth}.\qed

Now we investigate the case when the function $g$ in Corollary \ref{coro_ergodic} behaves as $z$ for small $z$ in the case of pure jump scaling. This type of functions make a highly oscillating drift term appear due to the small jumps. The fluctuations of that drift should overscale the size of the large jumps. However, when $g$ is odd, the fluctuations are stochastically centered (mean 0 w.r.t. $\mu$) so that we can establish the asymptotic convergence of these fluctuations towards their mean all the same:


\begin{theorem}{\bf Ergodic theorem II.}{\bf (Case of pure jump scaling).}\label{ergth2}
Consider a truncation function $h:\R\to \R$ such that $h(z)=z$ if $|z|\leq 1$ and $h(z)={\rm sign}(z)$ if $|z|>1$. Define ${\bfm h}_\epsilon\in L^\infty(\Omega)$ by 
$$ {\bfm h}_\epsilon(\omega)=\lim_{\alpha\downarrow 0}\frac{1}{\epsilon}\int_{|z|>\alpha}h(\delta(\epsilon)z){\bfm c}(\omega,z)e^{2{\bfm V}(\omega)}\chi(dz).$$
(To see why the limit exists, cf Lemma \ref{limitz}). Then
\begin{equation*}
\lim_{\epsilon\to 0}\M_\pi\E\Big[\big|\epsilon\int_0^{t/\epsilon}{\bfm h}_\epsilon(\tau_{X_{r-}}\omega)\,dr\big|\Big]=0.
\end{equation*}
\end{theorem}

\noindent {\bf Proof.} Choose a decreasing strictly positive sequence $(\beta_n)_{n\in\nat^*}$ converging towards $0$ as $n$ goes to $\infty$. For each $n\in \nat^*$, we define $h_n:\R\to \R$ as the truncation of $h$ at threshold $\beta_n$, that is  $h_n(z)=h(z)$ if $\beta_{n}<|z|$ and $0$ otherwise. Notice that $(h_n)_n$ uniformly converges towards $h$ on $\R$. We further define for each $n\in\nat$, $\alpha>0 $ and $\epsilon>0$,
\begin{align*}
 {\bfm h}^{n,\alpha}_\epsilon=&\frac{1}{\epsilon}\int_{|z|>\alpha}h_n(\delta(\epsilon)z){\bfm c}(\omega,z)e^{2{\bfm V}(\omega)}\chi(dz),\\
  {\bfm h}^{n}_\epsilon=&\frac{1}{\epsilon}\int_{\R}h_n(\delta(\epsilon)z){\bfm c}(\omega,z)e^{2{\bfm V}(\omega)}\chi(dz),\\
  {\bfm h}^{\alpha}_\epsilon=&\frac{1}{\epsilon}\int_{|z|>\alpha}h(\delta(\epsilon)z){\bfm c}(\omega,z)e^{2{\bfm V}(\omega)}\chi(dz)
 \end{align*} 
 The truncation w.r.t. $\alpha$ avoids dealing with integrability issues around $z=0$. Our strategy is the following. From Lemma \ref{implemma}, we can find a constant $C(n)$, only depending on $n$ and satisfying $\lim_{n\to \infty}C(n)=0$, such that 
 $$\M_\pi\E\Big|\epsilon\int_0^{t/\epsilon}({\bfm h}^\alpha_\epsilon-{\bfm h}^{n,\alpha}_\epsilon)(\tau_{X_{r-}}\omega)\,dr\Big|\leq C(n).$$
 Thus, Fatou's lemma yields
 \begin{align*}
\M_\pi\E\Big[\big|\epsilon\int_0^{t/\epsilon}{\bfm h}_\epsilon(\tau_{X_{r-}}\omega)\,dr\big|\Big]\leq &  \liminf_{\alpha\downarrow 0}\M_\pi\E\Big[\big|\epsilon\int_0^{t/\epsilon}{\bfm h}^\alpha_\epsilon(\tau_{X_{r-}}\omega)\,dr\big|\Big]\\ \leq &\liminf_{\alpha\downarrow 0}\M_\pi\E\Big[\big|\epsilon\int_0^{t/\epsilon}({\bfm h}^\alpha_\epsilon-{\bfm h}^{n,\alpha}_\epsilon)(\tau_{X_{r-}}\omega)\,dr\big|\Big]\\
&+\liminf_{\alpha\downarrow 0}\M_\pi\E\Big[\big|\epsilon\int_0^{t/\epsilon}{\bfm h}^{n,\alpha}_\epsilon(\tau_{X_{r-}}\omega)\,dr\big|\Big]\\ \leq &C(n)+\M_\pi\E\Big[\big|\epsilon\int_0^{t/\epsilon}{\bfm h}^{n}_\epsilon(\tau_{X_{r-}}\omega)\,dr\big|\Big]
 \end{align*}
 Clearly, we just have to prove that, for a fixed $n\in\nat^*$, $\M_\pi\E\Big[\big|\epsilon\int_0^{t/\epsilon}{\bfm h}^{n}_\epsilon(\tau_{X_{r-}}\omega)\,dr\big|\Big]\to 0$ as $\epsilon\to0$. This is a consequence of Corollary \ref{coro_ergodic} with ($g(z)=h_n(z)$). Indeed,  with $g(z)=h_n(z)$, the limit in  Corollary \ref{coro_ergodic} reduces to $0$, because the limit should match $t\int_{\R}h_n(z)\M[{\bfm \theta}(\cdot, z)]\mathcal{H}(dz)=\lim_{\epsilon\to 0}\frac{1}{\epsilon}\int_{\R}h_n(\delta(\epsilon)z)\M[{\bfm c}(\cdot,z)]\chi(dz)$. But the latter quantity is equal to $0$ since $h_n$ is odd, the measure $\chi$ is symmetric ($\chi(dz)=\chi(-dz)$) and $\M[{\bfm c}(\cdot,z)]$ is even by symmetry of ${\bfm c}$ (we have $\M[{\bfm c}(\cdot,-z)]=\M[{\bfm c}(\tau_z\cdot,-z)]=\M[{\bfm c}(\cdot,z)]$, $\chi$-a.s.). \qed

\begin{lemma}\label{implemma}
For any $n\in\nat$, $\alpha>0$ and $\epsilon>0$, we have
$$\forall\epsilon>0, \quad \M_\pi\Big[\big|\epsilon\int_0^{t/\epsilon}({\bfm h}^{n,\alpha}_\epsilon-{\bfm h}^{\alpha}_\epsilon)(\tau_{X_{r-}}\omega)\,dr\big|\Big]\leq C(\epsilon,n)$$ where $C(\epsilon,n)=\frac{\sup_{\Omega\times \R}|{\bfm c}|}{2}\frac{\delta(\epsilon)^2}{\epsilon}\int_{|z|\delta(\epsilon)\leq \beta_n}z^2\chi(dz)$. Moreover, from Assumption \ref{fouriersc}, we have $\lim_{n\to \infty}\sup_\epsilon C(\epsilon,n)=0$.
\end{lemma}

\noindent {\bf Proof.} We split the proof into 3 steps.

\noindent $\bullet$ Step 1: For $\epsilon>0$, $n\in \nat^*$, we define ${\bfm g}^{n,\alpha}_\epsilon={\bfm h}^{\alpha}_\epsilon -{\bfm h}^{n,\alpha}_\epsilon$. We claim: 
\begin{equation}\label{poincaren} 
\forall {\bfm \varphi}\in \mathcal{C},\,\,\,({\bfm g}_\epsilon^{n,\alpha},{\bfm \varphi})_\pi\leq  \big(\epsilon^{-1}C(\epsilon,n)\big)^{1/2}B^s({\bfm \varphi},{\bfm \varphi})^{1/2}.
\end{equation}

\noindent Proof. Since $h-h_{n}$ is odd, we use Lemma \ref{ippj} (with $g(z)=\one_{|z|>\alpha}(h-h_n)(\delta(\epsilon)z)$):
\begin{align*}
({\bfm g}^{n,\alpha}_\epsilon,{\bfm \varphi})_\pi&=\frac{1}{2\epsilon}\int_{|z|>\alpha}(h-h_{n})(z\delta(\epsilon))\M\big[{\bfm c}(\cdot,z)({\bfm \varphi}-T_{z}{\bfm \varphi})\big]\chi(dz)\\
&\leq \frac{1}{\sqrt{2}\epsilon}\Big(\M\int_{|z|>\alpha}(h-h_{n})^2(z\delta(\epsilon)){\bfm c}(\cdot,z)\chi(dz)\Big)^{1/2}B^j({\bfm \varphi},{\bfm \varphi})^{1/2}\\
&\leq \frac{1}{\epsilon}\Big(\frac{\sup_{\Omega\times \R}|{\bfm c}|}{2}\int_{\R}(h-h_{n})^2(z\delta(\epsilon))\chi(dz)\Big)^{1/2}B^s({\bfm \varphi},{\bfm \varphi})^{1/2}
\end{align*}
To conclude, it suffices to notice that $(h-h_{n})^2(z\delta(\epsilon))$ coincides with $\delta(\epsilon)^2z^2\one_{\{\delta(\epsilon)z\leq \beta_n\}}$ as soon as $\beta_n\leq 1$.

\noindent $\bullet$ Step 2: For each $n\in\nat $ and $\epsilon>0$, we define ${\bfm u}^{n,\alpha}_\epsilon=G_\epsilon({\bfm g}^{n,\alpha}_\epsilon)$. We claim:
\begin{equation}\label{bounds}
\epsilon^2|{\bfm u}^{n,\alpha}_\epsilon|_2^2+\epsilon({\bfm a}D{\bfm u}^{n,\alpha}_\epsilon,D{\bfm u}^{n,\alpha}_\epsilon)_\pi+\epsilon\M\int_{\R}|T_z{\bfm u}_\epsilon^{n,\alpha}-{\bfm u}_\epsilon^{n,\alpha}|^2{\bfm c}(\cdot;z)\chi(dz)\leq C(\epsilon,n,m)^2. 
\end{equation}

Proof. To see this, we just have to plug ${\bfm \psi}={\bfm u}^{n,\alpha}_\epsilon$ in the resolvent equation \eqref{eqbase} associated to ${\bfm g}^{n,\alpha}_\epsilon$. The right-hand side matches $({\bfm g}^{n,\alpha}_\epsilon,{\bfm u}^{n,\alpha}_\epsilon)_\pi$ and can be estimated as (see \eqref{poincaren})
$$({\bfm g}^{n,\alpha}_\epsilon,{\bfm u}^{n,\alpha}_\epsilon)_\pi\leq \big(\epsilon^{-1}C(\epsilon,n)\big)^{1/2}B^s(({\bfm u}^{n,\alpha}_\epsilon,{\bfm u}^{n,\alpha}_\epsilon)^{1/2}\leq \frac{C(\epsilon,n)^2}{2\epsilon}+\frac{B^s({\bfm u}^{n,\alpha}_\epsilon,{\bfm u}^{n,\alpha}_\epsilon)}{2}$$
so that the result follows by multiplying both sides of \eqref{eqbase} by $\epsilon$.

\noindent $\bullet$ Step 3: Since ${\bfm g}^{n,\alpha}_\epsilon\in H^\infty(\Omega)$, we have ${\bfm u}^{n,\alpha}_\epsilon\in H^\infty(\Omega) $ (cf Prop \ref{core}). Thus we apply the It\^o formula to the function ${\bfm u}^{n,\alpha}_\epsilon$ (cf Section \ref{particle}) and we get 
\begin{align*}
{\bfm u}^{n,\alpha}_\epsilon(\tau_{X_{t}}\omega) =& {\bfm u}^{n,\alpha}_\epsilon(\omega) +\int_0^t{\bfm L}{\bfm u}^{n,\alpha}_\epsilon(\tau_{X_{r-}}\omega)\,dr+\int_0^t{\bfm \sigma}D{\bfm u}^{n,\alpha}_\epsilon(\tau_{X_{r-}}\omega)\,dB_r\\&+\int_0^t\big({\bfm u}_\epsilon^{n,\alpha}(\tau_{X_{r-}+{\bfm \gamma(\tau_{X_{r-}}\omega,z)}}\omega)-{\bfm u}^{n,\alpha}_\epsilon(\tau_{X_{r-}}\omega)\big)\,d\tilde{N}(dr,dz)\\
=&{\bfm u}_\epsilon^{n,\alpha}(\omega)+\int_0^t(\epsilon{\bfm u}^{n,\alpha}_\epsilon-{\bfm h}^{n,\alpha}_\epsilon)(\tau_{X_{r-}}\omega)\,dr+\int_0^t{\bfm \sigma}D{\bfm u}^{n,\alpha}_\epsilon(\tau_{X_{r-}}\omega)\,dB_r\\&+\int_0^t\int_{\R}\big({\bfm u}_\epsilon^{n,\alpha}(\tau_{X_{r-}+{\bfm \gamma(\tau_{X_{r-}}\omega,z)}}\omega)-{\bfm u}_\epsilon^{n,\alpha}(\tau_{X_{r-}}\omega)\big)\,d\tilde{N}(dr,dz)
\end{align*}
We replace $t$ by $t/\epsilon$, multiply both sides of the above equality by $\epsilon$ and isolate the term corresponding to ${\bfm h}_\epsilon$. We get
\begin{align*}
\epsilon\int_0^{t/\epsilon}{\bfm g}^{n,\alpha}_\epsilon(\tau_{X_{r-}}\omega)\,dr=& \epsilon {\bfm u}^{n,\alpha}_\epsilon(\omega)-\epsilon{\bfm u}^{n,\alpha}_\epsilon(\tau_{X_{t/\epsilon}}\omega)+\epsilon^2\int_0^{t/\epsilon}{\bfm u}_\epsilon^{n,\alpha}(\tau_{X_{r-}}\omega)\,dr\\&+\epsilon\int_0^{t/\epsilon}{\bfm \sigma}D{\bfm u}^{n,\alpha}_\epsilon(\tau_{X_{r-}}\omega)\,dB_r\\
&+\epsilon\int_0^{t/\epsilon}\int_{\R}\big({\bfm u}_\epsilon^{n,\alpha}(\tau_{X_{r-}+{\bfm \gamma(\tau_{X_{r-}}\omega,z)}}\omega)-{\bfm u}_\epsilon^{n,\alpha}(\tau_{X_{r-}}\omega)\big)\,d\tilde{N}(dr,dz).
\end{align*}
The remaining part of the proof consists in proving that the quadratic mean of each term in the right-hand side of the above expression is bounded by  $C(\epsilon,n)^2$. The procedure is the same for each term: integrate the square of the term, use the invariance of $\pi$ for the process $Y(\omega)=\tau_X\omega$ and deduce the result from \eqref{bounds}. So we only detail the procedure for one term, say the last one.
\begin{align*}
\M_\pi\E\Big[\big|&\epsilon\int_0^{t/\epsilon}\int_{\R}\big({\bfm u}_\epsilon^{n,\alpha}(\tau_{X_{r-}+{\bfm \gamma(\tau_{X_{r-}}\omega,z)}}\omega)-{\bfm u}_\epsilon^{n,\alpha}(\tau_{X_{r-}}\omega)\big)\,d\tilde{N}(dr,dz)\big|^2\Big]\\
&\leq \M_\pi\E\Big[\epsilon^2\int_0^{t/\epsilon}\int_{\R}\big({\bfm u}_\epsilon^{n,\alpha}(\tau_{X_{r-}+{\bfm \gamma(\tau_{X_{r-}}\omega,z)}}\omega)-{\bfm u}^{n,\alpha}_\epsilon(\tau_{X_{r-}}\omega)\big)^2\nu(dz)dr\Big]\\
&= \epsilon^2\int_0^{t/\epsilon}\int_{\R}\M_\pi\E\Big[\big({\bfm u}_\epsilon^{n,\alpha}(\tau_{X_{r-}+{\bfm \gamma(\tau_{X_{r-}}\omega,z)}}\omega)-{\bfm u}^{n,\alpha}_\epsilon(\tau_{X_{r-}}\omega)\big)^2\Big]\nu(dz)dr\\
&= \epsilon t\M\Big[\int_{\R}(T_z{\bfm u}^{n,\alpha}_\epsilon-{\bfm u}^{n,\alpha}_\epsilon)^2{\bfm c}(\cdot,z)\chi(dz)\Big]\leq C(\epsilon,n)^2.\qed
\end{align*}
\section{Construction of the correctors}\label{sec:correctors}
In this section, we define the so-called correctors:

{\bf 1) Case of pure jump scaling.} No correctors. Actually, the job is already carried out in the proof of Th. \ref{ergth2}.

 {\bf 2) Case of diffusive scaling.} We define ${\bfm h}(\omega)=\lim_{\alpha\downarrow 0}\int_{|z|>\alpha}z{\bfm c}(\omega,z)\chi(dz)$  (Lemma \ref{limitz} together with $\int_\R z^2\chi(dz)<+\infty$ ensures the existence of the limit). Given $\lambda>0$,  we define $${\bfm u}_\lambda =G_\lambda({\bfm b}+{\bfm h})$$

\begin{remark}\label{regulcorr}
Since ${\bfm b}\in H^\infty(\Omega)$, $G_\lambda({\bfm b})\in H^\infty(\Omega)$ (see Proposition \ref{core}). Furthermore, from Lemma \ref{limitz} and the regularity conditions on ${\bfm c}$ (see Assumption \ref{regul}), it is plain to deduce that  ${\bfm h}\in H^\infty(\Omega)$ and the successive derivatives of ${\bfm h}$ are given, for $k\geq 1$, $D^k{\bfm h}=\lim_{\alpha\downarrow 0}\int_{|z|>\alpha}zD^k{\bfm c}(\omega,z)\chi(dz)$. Hence $G_\lambda({\bfm h})\in H^\infty(\Omega)$.
\end{remark}

\begin{proposition}\label{cvcorrector}{\bf Case of diffusive scaling.} There are $ {\bfm \xi}\in L^2(\Omega)$ and ${\bfm \zeta}\in L^2(\R\times\Omega;{\bfm c}(\omega;z)\chi(dz)d\mu(\omega))$ such that
$$\lambda|{\bfm u}_\lambda|^2_\pi+|D{\bfm u}_\lambda- {\bfm \xi}|^2_\pi+\M\int_{\R}|T_z{\bfm u}_\lambda-{\bfm u}_\lambda-{\bfm \zeta}(\cdot,z)|^2{\bfm c}(\cdot,z)\chi(dz)\to 0\quad \text{ as }\lambda\to 0.$$
\end{proposition}

\noindent {\bf Proof.} Remind that $\int_\R z^2 \chi(dz)<+\infty$. Applying Lemma \ref{ippj} (with $g(z)=z\one_{|z|>\alpha}$) and the Cauchy-Schwarz inequality yields: $\forall {\bfm v}\in \mathcal{C}$
$$({\bfm h},{\bfm v})_\pi=-\frac{1}{2}\lim_{\alpha\downarrow 0}\M\int_{|z|>\alpha}z{\bfm c}(\cdot,z)(T_z{\bfm v}-{\bfm v})\chi(dz)\leq \Big(\frac{\sup_{\Omega\times\R}|{\bfm c}|}{2}\int_\R z^2 \chi(dz)\Big)^{1/2}B^j({\bfm v},{\bfm v})^{1/2} .$$
By using integration by parts, we also get:
\begin{equation*}
({\bfm b},{\bfm v})_\pi=(1/2)(e^{2{\bfm V}}D(e^{-2{\bfm V}}{\bfm a}),{\bfm v})_\pi=-(1/2)({\bfm a},D{\bfm v})_\pi \leq \big(\M_\pi[{\bfm a}]/2\big)^{1/2}B^d({\bfm v},{\bfm v})^{1/2}
\end{equation*}
By gathering the above inequalities, we can find a constant $C$ such that
\begin{equation}\label{poincorr}
({\bfm b}+{\bfm h},{\bfm v})_\pi\leq CB^s({\bfm v},{\bfm v})^{1/2}.
\end{equation}
The standard inequality $ab\leq a^2/2+b^2/2 $ yields $({\bfm b}+{\bfm h},{\bfm v})_\pi\leq C^2/2+B^s({\bfm v},{\bfm v})/2$.
Plugging this in the right-hand side of (\ref{eqbase}), we get:
\begin{align*}
& \; \lambda({\bfm u}_\lambda,{\bfm v})_\pi+B^s({\bfm u}_\lambda,{\bfm v}) \leq C^2/2+B^s({\bfm v},{\bfm v})/2,
\end{align*}
from which one easily gets by setting ${\bfm v}={\bfm u}_\lambda$: $\lambda|{\bfm u}_\lambda|^2_\pi+B^s({\bfm u}_\lambda,{\bfm u}_\lambda) \leq C^2$. This implies the existence of $ {\bfm \xi}\in L^2(\Omega)$, ${\bfm \zeta}\in L^2(\R\times\Omega;{\bfm c}(\omega,z)\chi(dz)  d\mu(\omega))$ such that the following weak convergence holds along some subsequence:
\begin{equation}\label{convfaible}
D{\bfm u}_\lambda \underset{\lambda \to 0}{\rightarrow} {\bfm \xi}, \quad T_z{\bfm u}_\lambda-{\bfm u}_\lambda \underset{\lambda \to 0}{\rightarrow} {\bfm \zeta}.
\end{equation}
Actually, the convergence holds along the whole  subsequence since the limit is characterized by
\begin{equation}\label{weaklimit}
 \forall v\in \H,\quad ({\bfm b}+{\bfm h},{\bfm v})_\pi  
= \frac{1}{2}({\bfm a}{\bfm \xi},D{\bfm v})_\pi + \frac{1}{2}\M\int_{\R}{\bfm \zeta}(z,\cdot)(T_z{\bfm v}-{\bfm v}){\bfm c}(\cdot,z)\chi(dz),
\end{equation} which is obtained by letting $\lambda$ go to zero (along the subsequence) in \eqref{eqbase} (notice that $\lambda {\bfm u}_\lambda\to 0$ since $\lambda|{\bfm u}_\lambda|^2_\pi\leq C^2$).
By setting ${\bfm v}={\bfm u}_\delta$ above and letting $\delta$ go to zero, we have:
\begin{align*}
  \overline{\underset{\delta \to 0}{\lim}}({\bfm b},{\bfm u}_\delta)_\pi 
 & \leq \frac{1}{2}({\bfm a}{\bfm \xi},{\bfm \xi})_\pi + \frac{1}{2}\M\int_{\R}{\bfm \zeta}(z,\cdot)^2{\bfm c}(\cdot,z)\chi(dz)     
\end{align*}
Using once again relation (\ref{eqbase}) with ${\bfm \psi}={\bfm u}_\delta$, we conclude that:
\begin{align*}
& \overline{\underset{\delta \to 0}{\lim}} \left( \delta({\bfm u}_\delta,{\bfm u}_\delta)_\pi+\frac{1}{2}({\bfm a}D{\bfm u}_\delta,D{\bfm u}_\delta)_\pi+\frac{1}{2}\M\int_{\R}(T_z{\bfm u}_\delta-{\bfm u}_\delta)^2{\bfm c}(\cdot,z)\chi(dz)   \right)  \\
& \leq \frac{1}{2}({\bfm a}{\bfm \xi},{\bfm \xi})_\pi + \frac{1}{2}\M\int_{\R}{\bfm \zeta}(z,\cdot)^2{\bfm c}(\cdot,z)\chi(dz) .
\end{align*}
From this, we deduce that the weak convergences in (\ref{convfaible}) are in fact strong and that $\overline{\underset{\delta \to 0}{\lim}} \; \delta |{\bfm u}_\delta|^2_2=0$.\qed

\section{Tightness}\label{tightness}

Our strategy to establish the tightness of the "environment as seen from the particle" does not differ from \cite[Section 3.3]{olla} (idea originally due to \cite{varadhan}) and relies on the so-called Garcia-Rodemich-Rumsey inequality. So we set out the main steps of the proof, only proving what differs from \cite{olla} (only minor things), and let the reader be referred to \cite{olla} for further details. 

\begin{remark}
The setup in \cite{olla} is more general than ours in the sense that the author considers possibly non-symmetric processes. To simplify the reading, take ${\bfm A}=0$ in \cite{olla}. 
\end{remark}

More precisely, our pupose is the following
\begin{theorem}
Consider a family of functions $({\bfm h}_\epsilon)_\epsilon\subset L^\infty(\Omega)$ satisfying the following estimate:
\begin{equation}\label{poincaretight}
\forall {\bfm \varphi}\in \mathcal{C},\quad ({\bfm h}_\epsilon,{\bfm \varphi}^2)_\pi\leq P B^s( {\bfm \varphi}, {\bfm \varphi})^{1/2}| {\bfm \varphi}|_\pi
\end{equation}
for some positive constant $P$. Then we can establish the following continuity modulus estimate:
\begin{equation}
\M_\pi\E\Big[\sup_{\substack{|t-s|\leq \delta\\0\leq s,t\leq T}}\big|\epsilon\int_{s/\epsilon^2}^{t/\epsilon^2}{\bfm h}_\epsilon(Y_{r-}(\omega))\,dr\big|\Big]\leq C(T)\delta^{1/2}\ln \delta^{-1}
\end{equation}
for some positive constant $C(T)$ only depending on $T$.
\end{theorem}

\noindent {\bf Guideline of the proof.} To begin with, we remind the reader of the GRR inequality:
\begin{proposition}\label{rumsey}{\bf (Garsia-Rodemich-Rumsey's
inequality).} Let $p$ and $\Psi$ be strictly increasing continuous
functions on $[0,+\infty[$ satisfying $p(0)=\Psi(0)=0$ and
$\lim_{t\to \infty}\Psi(t)=+\infty$. For given $T>0$ and $g\in
C([0,T];\R^d)$, suppose that there exists a finite $B$ such that;
\begin{equation}\label{GRRcond}
  \int_0^T\int_0^T\Psi\Big(\frac{|g(t)-g(s)|}{p(|t-s|)}\Big)\,ds\,dt\leq B<\infty.
\end{equation}
Then, for all $0\leq s \leq t \leq T$,
\begin{equation}\label{GRR}
  |g(t)-g(s)|\leq 8\int_0^{t-s}\Psi^{-1}(4B/u^2)\,dp(u).
\end{equation}
\end{proposition}

The first step is to estimate the exponential moments of the random variable $\epsilon\int_{s/\epsilon^2}^{t/\epsilon^2}{\bfm h}_\epsilon(Y_{r-}(\omega))\,dr$. It turns out that the Feynmann-Kac formula provides a connection between the exponential moments and the solution of a certain evolution equation:
\begin{theorem}{\bf Feynmann-Kac formula.}\label{fk}
Let ${\bfm U}$ belong to $L^\infty(\Omega)$. Then the function
$${\bfm u}(t,\omega)=\E\Big[\exp\big(\int_0^t{\bfm U}(Y_{r-}(\omega))\,dr\big)\Big]$$ is a solution of the equation
$$ \partial_t{\bfm u}={\bfm L}{\bfm u}+{\bfm U}{\bfm u}$$ with initial condition ${\bfm u}(0,\omega)={\bfm 1}$. 
\end{theorem}
\begin{remark}
By solution, we mean a function ${\bfm u}$ such that $\forall t \geq 0$, ${\bfm u}(t,\cdot)\in{\rm Dom}({\bfm L})$ and $$\lim_{s\to 0}\frac{{\bfm u}(t+s,\cdot)-{\bfm u}(t,\cdot)}{s}={\bfm L}{\bfm u}(t,\cdot)+{\bfm U}(\cdot){\bfm u}(t,\cdot)\quad \text{in }L^2(\Omega).$$
\end{remark}
\begin{remark} Though it is not necessary, the author also proves in \cite[Theorem 3.2]{olla} uniqueness of the solution to the equation. So, the reader may skip the corresponding part of the proof.
\end{remark}
Using the equation satisfied by ${\bfm u}(t,\cdot)$, we are now in position to establish bounds for the function ${\bfm u}$
\begin{proposition}\label{carreu}
Let ${\bfm u}(t,\cdot)$ be the function of Theorem \ref{fk}. Then
$$ \M_\pi[{\bfm u}(t,\cdot)^2]\leq e^{2\lambda_0({\bfm L}+{\bfm U})t}$$ where $\lambda_0({\bfm L}+{\bfm U})$ is defined as  $\lambda_0({\bfm L}+{\bfm U})=\sup_{\substack{|{\bfm \varphi}|_\pi=1,\\{\bfm \varphi}\in{\rm Dom}{\bfm L}}}({\bfm \varphi},({\bfm L}+{\bfm U}){\bfm \varphi})_\pi$.
\end{proposition} 
 Following \cite[Theorem 3.4]{olla}, we make use of Proposition \ref{carreu} to prove
 $$\M_\pi\E\Big[\exp\Big|\alpha\epsilon\int_{s/\epsilon^2}^{t/\epsilon^2}{\bfm U}(Y_{r-}(\omega))\,dr\Big|\Big] \leq 2\exp\Big(\lambda_0(\epsilon^{-2}{\bfm L}+\epsilon^{-1}\alpha{\bfm U})(t-s)\Big).$$
 In particular, we can choose ${\bfm U}={\bfm h}_\epsilon$ and use \eqref{poincaretight} to get $\lambda_0(\epsilon^{-2}{\bfm L}+\epsilon^{-1}\alpha{\bfm h}_\epsilon)\leq \alpha^2P^2/4$. This yields
$$\M_\pi\E\Big[\exp\Big|\alpha\epsilon\int_{s/\epsilon^2}^{t/\epsilon^2}{\bfm h}_\epsilon(Y_{r-}(\omega))\,dr\Big|\Big] \leq 2 \exp\Big(\alpha^2P^2(t-s)/4\Big).$$
We conclude by using the GRR inequality (with $g(t)=\epsilon\int_0^{t/\epsilon^2}{\bfm h}_\epsilon(Y_{r-}(\omega))\,dr$, $p(t)=\sqrt{t}$, $\Psi(t)=e^t-1$), by taking the expectation and by using the above estimate.\qed

We conclude this section by making three important remarks. First, notice that ${\bfm b}$ satisfies the relation \eqref{poincaretight} since for any ${\bfm \varphi}\in \mathcal{C}$
$$({\bfm b},{\bfm \varphi}^2)_\pi=(D(e^{-2{\bfm V}}{\bfm a}),{\bfm \varphi}^2)_2=-2({\bfm a},{\bfm \varphi}D{\bfm \varphi})_\pi\leq 2|{\bfm a}|_\infty^{1/2}({\bfm a}D{\bfm \varphi},D{\bfm \varphi})_\pi^{1/2}|{\bfm \varphi}|_\pi.$$
We deduce
\begin{equation}\label{tightnessb}
\M_\pi\E\Big[\sup_{\substack{|t-s|\leq \delta\\0\leq s,t\leq T}}\big|\epsilon^{1/2}\int_{s/\epsilon}^{t/\epsilon}{\bfm b}(Y_{r-}(\omega))\,dr\big|\Big]\leq C(T)\delta^{1/2}\ln \delta^{-1}.
\end{equation}
Second, define the function $h$ by $h(z)=z$ if $|z|\leq 1$, $h(z)={\rm sign}(z)$ if $|z|>1$, and (the limit exists in the $L^\infty$ sense because of Lemma \ref{limitz}) 
\begin{equation}\label{he}
{\bfm h}_\epsilon=\lim_{\alpha\downarrow 0}\frac{1}{\epsilon}\int_{|z|>\alpha}h(z\delta(\epsilon)){\bfm c}(\cdot,z)e^{2{\bfm V}}\chi(dz).
\end{equation}
Since $h$ is odd, we can apply Lemma \ref{ippj} to obtain: for any ${\bfm \varphi}\in \mathcal{C}$
\begin{align*}
\epsilon^{\frac{1}{2}}({\bfm h}_\epsilon,{\bfm \varphi}^2)_\pi&=\lim_{\alpha\downarrow 0}\frac{1}{2\epsilon^{\frac{1}{2}}}\M\int_{|z|>\alpha}h(z\delta(\epsilon)){\bfm c}(\cdot,z){\bfm \varphi}^2\chi(dz)\\
&=-\lim_{\alpha\downarrow 0}\frac{1}{2\epsilon^{\frac{1}{2}}}\int_{|z|>\alpha}h(z\delta(\epsilon)){\bfm c}(\cdot,z)(T_{z}{\bfm \varphi}^2-{\bfm \varphi}^2)\chi(dz)\\
&\leq (2\epsilon)^{-\frac{1}{2}}\Big(\M\int_{|z|>0}h^2(z\delta(\epsilon)){\bfm c}(\cdot,z)({\bfm \varphi}+T_{z}{\bfm \varphi})^2\chi(dz)\Big)^{1/2}B^j({\bfm \varphi},{\bfm \varphi})^{1/2}\\
&\leq \big(2\sup_{\Omega\times \R}|{\bfm c}|\big)^{1/2}\epsilon^{-\frac{1}{2}}\Big(\int_{\R}h^2(z\delta(\epsilon))\chi(dz)\Big)^{1/2}B^s({\bfm \varphi},{\bfm \varphi})^{1/2}|{\bfm \varphi}|_\pi.
\end{align*}
In the case of pure jump scaling, the quantity $\epsilon^{-1/2}\big(\int_{\R}h^2(z\delta(\epsilon))\chi(dz)\big)^{1/2}$ is bounded by a constant  independent of $\epsilon$ (see Assumption \ref{fouriersc}). So, we can apply our estimates to the function $\epsilon^{\frac{1}{2}}{\bfm h}_\epsilon$ and get
\begin{equation}\label{tightnessh}
\M_\pi\E\Big[\sup_{\substack{|t-s|\leq \delta\\0\leq s,t\leq T}}\big|\epsilon\int_{s/\epsilon}^{t/\epsilon}{\bfm h}_\epsilon(Y_{r-}(\omega))\,dr\big|\Big]\leq C(T)\delta^{1/2}\ln \delta^{-1}.
\end{equation}
Third, in the case of diffusive scaling, that is $\int_{\R}z^2\chi(dz)<+\infty$, we consider the function 
${\bfm h}=\lim_{\alpha\downarrow 0}\int_{|z|>\alpha}z{\bfm c}(\cdot,z)e^{2{\bfm V}}\chi(dz)$ (see Lemma \ref{limitz} again concerning the existence of the limit). Once again, by applying Lemma \ref{ippj}, we can derive the following estimate:
$$({\bfm h},{\bfm \varphi}^2)_\pi\leq \Big(\frac{1}{2}\M\int_\R(T_z{\bfm \varphi}+{\bfm \varphi})^2z^2{\bfm c}(\cdot,z)\chi(dz)\Big)^{1/2}B^j({\bfm \varphi},{\bfm \varphi})^{1/2}\leq (2\sup|{\bfm c}|)^{1/2}|{\bfm \varphi}|_\pi B^j({\bfm \varphi},{\bfm \varphi})^{1/2},$$ from which we deuce
\begin{equation}\label{tightnesshdiff}
\M_\pi\E\Big[\sup_{\substack{|t-s|\leq \delta\\0\leq s,t\leq T}}\big|\epsilon^{1/2}\int_{s/\epsilon}^{t/\epsilon}{\bfm h}(Y_{r-}(\omega))\,dr\big|\Big]\leq C(T)\delta^{1/2}\ln \delta^{-1}.
\end{equation}

\section{Homogenization}


In this section, we prove the homogenization theorem. 

{\bf 1) Case of pure jump scaling.}  From \eqref{SDE}, we have the following equation for the rescaled process $\delta(\epsilon)X_{\cdot/\epsilon}$:
\begin{align*}
\delta(\epsilon)X_{t/\epsilon}=&\delta(\epsilon)\int_0^{t/\epsilon}  {\bfm b}(\tau_{X_{r-}}\omega) dr+\delta(\epsilon)\int_0^{t/\epsilon}  {\bfm e}(\tau_{X_{r-}}\omega) dr+\delta(\epsilon)\int_0^{t/\epsilon}{\bfm \sigma}(\tau_{X_{r-}}\omega)\,dB_r\\
&+\delta(\epsilon)\int_0^{t/\epsilon}\int_{\R}\gamma(\tau_{X_{r-}}\omega,z)\,\hat{N}(dr,dz) . 
\end{align*}

In order to prove the result, we consider each term in the above sum separately. In view of \eqref{tightnessb}, we have
$$\M_\pi\E\Big[\sup_{0\leq t\leq T}\big|\delta(\epsilon)\int_{0}^{t/\epsilon}{\bfm b}(Y_{r-}(\omega))\,dr\big|\Big]\leq \frac{\delta(\epsilon)}{\epsilon^{1/2}}C(T)T^{1/2}\ln T^{-1} \to 0,\quad \text{as }\epsilon\to 0.$$
Concerning the Brownian martingale, by using the invariance of the measure $\pi$ for the process $Y(\omega)=\tau_X\omega $, we have
\begin{align*}
\M_\pi\E\Big[\sup_{0\leq t\leq T}\big|\delta(\epsilon)\int_0^{t/\epsilon}{\bfm \sigma}(\tau_{X_{r-}}\omega)\,dB_r\big|^2\Big]\leq  \M_\pi\E\Big[\delta(\epsilon)^2\int_0^{T/\epsilon}{\bfm a}(\tau_{X_{r-}}\omega)\,dr\Big]\leq \frac{\delta(\epsilon)^2}{\epsilon}T\M_\pi[{\bfm a}] .
\end{align*}
Thus, we just have to investigate the convergence of the following semimartingale $Y^{\epsilon}_t$:
\begin{equation*}
Y^{\epsilon}_t=\delta(\epsilon) \int_0^{ t/ \epsilon } {\bfm e}(\tau_{X_{r-}}\omega)dr+\delta(\epsilon) \int_0^{ t/ \epsilon } \int_{\R}\gamma(\tau_{X_{r-}}\omega,z)\,\hat{N}(dr,dz)
\end{equation*}
In order to obtain the desired result, we introduce the truncation function $h$ as defined in the Ergodic theorem \ref{ergth2} and we use theorem VIII.4.1 in \cite{jacod}. Following the notations of \cite{jacod}, we introduce the following processes:
\begin{equation*}
\check{Y}^{\epsilon,(h)}_t=\sum_{0<s \leq t}\Delta Y^{\epsilon}_s-h(\Delta Y^{\epsilon}_s)
\end{equation*}
and 
\begin{equation*}
Y^{\epsilon,(h)}_t=Y^{\epsilon}_t-\check{Y}^{\epsilon,(h)}_t.
\end{equation*}
Note that we can decompose the semimartingale $Y^{\epsilon,(h)}$ as:
\begin{equation*}
Y^{\epsilon}_t=M^{\epsilon,(h)}_t+B^{\epsilon,(h)}_t,
\end{equation*}
where $M^{\epsilon,(h)},B^{\epsilon,(h)}$ are given by:
\begin{equation*}
M^{\epsilon,(h)}_t= \int_0^{ t/ \epsilon } \int_{\R}  h(\delta(\epsilon)\gamma(\tau_{X_{r-}}\omega,z))\,\tilde{N}(dr,dz)
\end{equation*}
and 
\begin{align*}
B^{\epsilon,(h)}_t= &\delta(\epsilon)\int_0^{t/ \epsilon}{\bfm e}(\tau_{X_{r-}}\omega)dr+\int_0^{t/ \epsilon} \int_{|z|> 1} h\big(\delta(\epsilon){\bfm \gamma} (\tau_{X_{r-}}\omega,z)\big)\nu(dz)dr.
\end{align*}
As soon as $\delta(\epsilon)S\leq 1$ (cf Assumption \ref{regul}.4 for the definition of $S$), we have 
$$\delta(\epsilon){\bfm e}(\omega)=\lim_{\alpha\downarrow 0}\int_{\alpha\leq |{\bfm \gamma}|}\delta(\epsilon){\bfm \gamma}(\omega,z)\one_{\{|z|\leq 1\}}\nu(dz)=\lim_{\alpha\downarrow 0}\int_{\alpha\leq |{\bfm \gamma}|}h\big(\delta(\epsilon){\bfm \gamma}(\omega,z)\big)\one_{\{|z|\leq 1\}}\nu(dz),$$
in such a way that $B^{\epsilon,(h)}$ can be rewritten as (cf the notations of Theorem \ref{ergth2})
$$B^{\epsilon,(h)}_t= \epsilon\int_0^{t/\epsilon}{\bfm h}_\epsilon(\tau_{X_{r-}}\omega)dr.$$
According to \eqref{tightnessh}, $B^{\epsilon,(h)}$ is tight in $D(\R_+;\R)$ for the Skorohod topology. Moreover, Theorem \ref{ergth2} ensures that the finite-dimensional distributions of $B^{\epsilon,(h)}$ converges to $0$. Hence, $B^{\epsilon,(h)}$ converges to $0$ in probability in $D(\R_+;\R)$.

By Corollary \ref{coro_ergodic}, we have also the following convergence for $<M^{\epsilon,(h)}>_t$:
\begin{align*}
<M^{\epsilon,(h)}>_t & = \int_0^{t/ \epsilon} \int_{\R} h(\delta(\epsilon)z)^2{\bfm c}(\tau_{X_{r-}}\omega,z)e^{2{\bfm V}(\tau_{X_{r-}}\omega)}\chi(dz)dr    \\
& \overset{\epsilon \to 0+}{\longrightarrow} t\int_{\R}h(z)^2\M[{\bfm \theta}(\omega,{\rm sign}(z))]\mathcal{H}(dz)
\end{align*}
To sum up, the three characteristics of the semimartingale $Y^\epsilon$ converge as $\epsilon\to 0$ to those of a L\'evy process $L$ with L\'evy exponent:
\begin{equation*}
\varphi(u)=\int_{\R}(e^{iuz}-1-iuz\one_{\{|z|\leq 1\}})\M[{\bfm \theta}(\omega,{\rm sign}(z))]\mathcal{H}(dz).
\end{equation*}
Using theorem VIII.4.1 in \cite{jacod}, we conclude that the following convergence holds for the Skorohod topology:
\begin{equation*}
Y^{\epsilon} \overset{\epsilon \to 0}{\longrightarrow} L.\qed
\end{equation*}

\noindent {\bf Case of diffusive scaling}. We apply the It\^o formula to the function ${\bfm u}_\epsilon=G_\epsilon({\bfm b}+{\bfm h})$:
\begin{align}
{\bfm u}_\epsilon&(\tau_{X_t}\omega)-{\bfm u}_\epsilon(\omega)=\int_0^t \epsilon {\bfm u}_\epsilon(\tau_{X_{r-}}\omega) dr -\int_0^t {\bfm b}(\tau_{X_{r-}}\omega) dr -\int_0^t {\bfm h}(\tau_{X_{r-}}\omega) dr\nonumber\\
\label{itou}&+\int_0^t\int_\R\big({\bfm u}_\epsilon(\tau_{X_{r-}+{\bfm \gamma}(\tau_{X_{r-}}\omega,z)}\omega)-{\bfm u}_\epsilon (\tau_{X_{r-}}\omega)\big)\,\tilde{N}(dr,dz)+\int_0^tD{\bfm u}_\epsilon{\bfm \sigma}(\tau_{X_{r-}}\omega)\,dB_r. 
\end{align}
Therefore, by summing with \eqref{SDE} and by using the relation $${\bfm e}(\omega)-{\bfm h}(\omega)=-\int_{|z|>1}{\bfm \gamma}(\omega,z)\nu(dz),$$ we deduce:
\begin{align*}
{\bfm u}_\epsilon&(\tau_{X_t}\omega)+X_t={\bfm u}_\epsilon(\omega)+\int_0^t \epsilon {\bfm u}_\epsilon(\tau_{X_{r-}}\omega) dr +\int_0^t(1+D{\bfm u}_\epsilon){\bfm \sigma}(\tau_{X_{r-}}\omega)\,dB_r\\
&+\int_0^t\int_\R\big(\gamma(\tau_{X_{r-}}\omega,z)+{\bfm u}_\epsilon(\tau_{X_{r-}+{\bfm \gamma}(\tau_{X_{r-}}\omega,z)}\omega)-{\bfm u}_\epsilon (\tau_{X_{r-}}\omega)\big)\,\tilde{N}(dr,dz). 
\end{align*}
We now analyze the convergence of each rescaled term of the above relation. By Prop \ref{cvcorrector}, we have:
\begin{equation}\label{tfac}
\M_\pi\E\big[\sup_{0\leq t \leq T}\big|\delta(\epsilon)\int_0^{t/\epsilon} \epsilon {\bfm u}_\epsilon(\tau_{X_{r-}}\omega) dr \big|\big]\leq \delta(\epsilon)|{\bfm u}_\epsilon|_\pi\to 0,\quad \epsilon\to 0. 
\end{equation}
We know focus on $\delta(\epsilon)\big({\bfm u}_\epsilon(\tau_{X_t}\omega)-{\bfm u}_\epsilon(\omega)\big)$. Prop. \ref{cvcorrector} leads to $$\M_\pi\E\big[\big|\delta(\epsilon)\big({\bfm u}_\epsilon(\tau_{X_t}\omega)-{\bfm u}_\epsilon(\omega) \big|\big]\leq 2\delta(\epsilon)|{\bfm u}_\epsilon|_\pi\to 0,\quad \epsilon\to 0.$$
To see why the process $\delta(\epsilon)\big({\bfm u}_\epsilon(\tau_{X_t}\omega)-{\bfm u}_\epsilon(\omega)\big)$ is tight for the Skorohod topology, we have to get back to \eqref{itou}. In the right-hand side, we have already establish the tightness of all the terms with bounded variations (cf \eqref{tfac} \eqref{tightnessb} and \eqref{tightnesshdiff}). Concerning the martingale terms, it suffices to apply Corollary \ref{coro_ergodic} together with Prop \ref{cvcorrector} to the brackets to show that they converge to a continuous deterministic process (for further details, see the argument below). Hence the martingale terms are also tight, and so is $\delta(\epsilon)\big({\bfm u}_\epsilon(\tau_{X_t}\omega)-{\bfm u}_\epsilon(\omega)\big)$. To sum up, it  converges in probability for the Skorohod topology towards $0$.

It remains to treat the martingale term
\begin{align*}
M^\epsilon_t=&\epsilon^{1/2}\int_0^{t/\epsilon}\int_\R\big(\gamma(\tau_{X_{r-}}\omega,z)+{\bfm u}_\epsilon(\tau_{X_{r-}+{\bfm \gamma}(\tau_{X_{r-}}\omega,z)}\omega)-{\bfm u}_\epsilon (\tau_{X_{r-}}\omega)\big)\,\tilde{N}(dr,dz)\\
&+\epsilon^{1/2}\int_0^{t/\epsilon}(1+D{\bfm u}_\epsilon){\bfm \sigma}(\tau_{X_{r-}}\omega)\,dB_r.
\end{align*}
By using Proposition \ref{cvcorrector} and Corollary \ref{coro_ergodic}, the brackets 
\begin{align*}
<M^\epsilon>_t=&\epsilon\int_0^{t/\epsilon}\int_\R\big(z+T_z{\bfm u}_\epsilon-{\bfm u}_\epsilon\big)^2 (\tau_{X_{r-}}\omega){\bfm c}(\tau_{X_{r-}}\omega,z)\chi(dz)+\epsilon\int_0^{t/\epsilon}(1+D{\bfm u}_\epsilon)^2{\bfm a}(\tau_{X_{r-}}\omega)\,dr
\end{align*}
converge to the continuous deterministic process $t\mapsto At$ ($A$ is given by \eqref{form:A}). Using the martingale central limit theorem, cf \cite{helland}, we see that $(M^\epsilon)_\epsilon$ converges in law towards a Brownian motion with covariance matrix $A$ (note that the jump condition required in \cite{helland} results from Corollary \ref{coro_ergodic}) .\qed

\appendix
\section*{Appendix}
\section{Auxiliary lemmas}

\begin{lemma}\label{ippj}
Let $g:\R\to\R$ be a $\chi$-integrable  odd function, and let ${\bfm h}$ be defined as ${\bfm h}(\omega)=\int_\R g(z){\bfm c}(\omega,z)e^{2{\bfm V}}\chi(dz)$. Then, for any ${\bfm \varphi}\in \mathcal{C}$
$$({\bfm h},{\bfm \varphi})_\pi=-\frac{1}{2}\M\int_\R g(z)(T_z{\bfm \varphi}-{\bfm \varphi}){\bfm c}(\cdot,z)\chi(dz). $$
\end{lemma}

\noindent {\bf Proof.} We have to use the symmetry of ${\bfm c}$ ($\chi(dz)$ a.s., $2{\bfm c}(\omega,z)= {\bfm c}(\tau_z\omega,-z)+{\bfm c}(\omega,z)$) and the symmetry of $\chi$ ($\chi(dz)=\chi(-dz)$):
\begin{align*}
({\bfm h},{\bfm \varphi})_\pi&=\frac{1}{2}\int_{\R}g(z)\M[T_z{\bfm c}(\cdot,-z)+{\bfm c}(\cdot,z)){\bfm \varphi}]\chi(dz)\\
&=\frac{1}{2}\int_{\R}g(z)\M[{\bfm c}(\cdot,-z),T_{-z}{\bfm \varphi}]\chi(dz)+\frac{1}{2}\int_{\R}g(z)\M[{\bfm c}(\cdot,z),{\bfm \varphi}]\chi(dz)\\
&=-\frac{1}{2}\int_{\R}g(z)\M[{\bfm c}(\cdot,z),T_{z}{\bfm \varphi}]\chi(dz)+\frac{1}{2}\int_{\R}g(z)\M[{\bfm c}(\cdot,z),{\bfm \varphi}]\chi(dz)\\
&=\frac{1}{2}\M\int_{\R}g(z){\bfm c}(\cdot,z)({\bfm \varphi}-T_{z}{\bfm \varphi})\chi(dz)\qed
\end{align*}

\begin{lemma}\label{limitz}
Fix $k\in\nat$. If a measurable function ${\bfm g}:\Omega\times\R\to \R$ satisfies $$|{\bfm g}(\omega,z)|\one_{\{|z|\leq 1\}}\leq {\bfm C}(\omega)|z|,\quad |{\bfm g}(\omega,z)+{\bfm g}(\omega,-z)|\one_{\{|z|\leq 1\}}\leq {\bfm C}(\omega)|z|^2$$ for some function ${\bfm C}\in L^2(\Omega)$ (resp. ${\bfm C}\in L^\infty(\Omega)$) then the following limit exists in the $L^2(\Omega)$-sense (resp. $L^\infty(\Omega)$-sense):
$$\lim_{\alpha\downarrow 0}\int_{\alpha<|z|\leq 1}{\bfm g}(\omega,z)D^k{\bfm c}(\omega,z)\chi(dz).$$
\end{lemma}

\noindent {\bf Proof.} First notice that $D^k{\bfm c}$ is symmetric because ${\bfm c}$ is, that is $D^k{\bfm c}(\tau_z\omega,-z)=D^k{\bfm c}(\omega,z)$ $\chi(dz)$ a.s. In particular, since the mapping $x\in\R\mapsto D^k{\bfm c}(\tau_x\omega,z)$ is smooth, we have $ \chi(dz)$ a.s. 
$$D^k{\bfm c}(\tau_z\omega,-z)=D^k{\bfm c}(\omega,-z)+z\int_0^1D^{k+1}{\bfm c}(\tau_{zu}\omega,-z)du.$$
By plugging this into the relation $D^k{\bfm c}(\omega,z)=\frac{1}{2} \big(D^k{\bfm c}(\tau_z\omega,-z)+D^k{\bfm c}(\omega,z)\big)$, it is plain to see that, for $\alpha>0$:
\begin{align*}
\int_{\alpha<|z|\leq 1}{\bfm g}(\omega,z)D^k{\bfm c}(\omega,z)\chi(dz)=&\frac{1}{2}\int_{\alpha<|z|\leq 1}{\bfm g}(\omega,z)\big(D^k{\bfm c}(\omega,-z)+D^k{\bfm c}(\omega,z)\big)\chi(dz)\\&+\frac{1}{2}\int_{\alpha<|z|\leq 1}z{\bfm g}(\omega,z)\int_0^1D^{k+1}{\bfm c}(\tau_{zu}\omega,-z)du\chi(dz)\\=&\frac{1}{2}\int_{\alpha<|z|\leq 1}D^k{\bfm c}(\omega,z)\big({\bfm g}(\omega,z)+{\bfm g}(\omega,-z)\big)\chi(dz)\\&+\frac{1}{2}\int_{\alpha<|z|\leq 1}z{\bfm g}(\omega,z)\int_0^1D^{k+1}{\bfm c}(\tau_{zu}\omega,-z)du\chi(dz).
\end{align*}
We complete the proof thanks to  the bounds $|D^k{\bfm c}(\cdot,z)|_\infty+|D^{k+1}{\bfm c}(\cdot,z)|_\infty\leq C_k+C_{k+1}$ (see ssumption \ref{regul}.2), the estimates on ${\bfm g}$ and the relation $\int_{|z|\leq 1}z^2\chi(dz)<+\infty$.\qed

\begin{lemma}\label{contdc}
Consider a kernel ${\bfm d}:\Omega\times\R\to \R$ such that there is a constant $M\geq 0$ satisfying $|{\bfm d}(\cdot,z)|_\infty\leq M$ $\chi(dz)$ a.s. For each ${\bfm \varphi},{\bfm \psi}\in\H$ we have
$$\M\int_{\R^*}(T_z{\bfm \varphi}-{\bfm \varphi})(T_z{\bfm \psi}-{\bfm \psi}){\bfm d}(\cdot,z)\chi(dz) \leq C\big(({\bfm \varphi},{\bfm \psi})_2+(D{\bfm \varphi},D{\bfm \psi})_2\big)$$ for some constant $C\geq 0$ only depending on $ M$ and $ \chi$. 
\end{lemma}

\noindent {\bf Proof.} It suffices to split the integral w.r.t. the variable $z$ in two parts: for $|z|\leq1$ and $|z|> 1$. The first integral is estimated with the derivative  $D{\bfm \varphi}$, whereas the second is estimated with ${\bfm \varphi}$. Since that type of result is quite classical, details are left to the reader.\qed

\section{Proofs of Section \ref{statement}}

\noindent {\bf Proof of Lemma \ref{constructiongamma}.} We split the proof into several steps:

$\bullet$ \textit{Construction of ${\bfm \gamma}$ and $\nu$}:
  We define ${\bfm h}(\omega,z)= \int_z^{+\infty} {\bfm c}^s(\omega,r)\chi(r)\,dr $ if $z>0$ and  ${\bfm h}(\omega,z)=-\int_{-\infty}^z {\bfm c}^s(\omega,r)\chi(r)\,dr$ if $z<0$. We also define $F(z)=M\int_z^{+\infty}\chi(r)\,dr$ if $z>0$ and $F(z)=-M\int_{-\infty}^z \chi(r)\,dr$ if $z<0$.  Notice that, for any fixed $\omega$, ${\bfm h}(\omega,\cdot)$ and $F$ are both homeomorphisms  from $\R^*_+$ onto itself and from $\R^*_-$ onto itself.

Set ${\bfm \gamma}(\omega,z)={\bfm h}^{-1}(\omega, F(z))$ for $z\not =0$, which can be continuously extended by setting $ {\bfm \gamma}(\omega,0)=0$, and $\nu(z)=M\chi(z)$ for $z\in\R$. We should point out that, for each fixed $\omega$, the mapping $z\mapsto {\bfm \gamma}(\omega,z)$ is a homeomorphism from $\R$ onto itself.

Fix $\omega\in\Omega$. For $z>0$, ${\bfm \gamma}(\omega,\cdot)$ satisfy the relation ${\bfm h}(\omega, z)=F({\bfm \gamma}^{-1}(\omega,z))$, that is $\int_z^{+\infty}{\bfm c}^s(\omega,r)\chi(r)dr=\nu\big({\bfm \gamma}^{-1}(\omega,\cdot)([z,+\infty[)\big)$. Since the sets $[z,+\infty[$ for $z>0$ generate the Borelian $\sigma$-field of $]0,+\infty[$, the measures $\nu\circ {\bfm \gamma}^{-1}(\omega,\cdot)$ and ${\bfm c}^s(\omega,z)\chi(r)dz $ coincide on $]0,+\infty[$. Similarly, we prove that they coincide on $]-\infty,0[$, hence on $\R$.

Furthermore, notice that ${\bfm \gamma}$ satisfies the relation $$F(z)=h(\omega, {\bfm \gamma}(\omega,z))\leq  F({\bfm \gamma}(\omega,z)).$$
Since $F$ is strictly decreasing on $\R_+^*$, we deduce $z\geq {\bfm \gamma}(\omega,z)$ for $z>0$. Since for $z>0$, ${\bfm \gamma}(\omega,z)>0$, we deduce $|{\bfm \gamma}(\omega,z)|\leq |z|$ for $z>0$.
The same estimate holds for $z<0$ in such a way that $|{\bfm \gamma}(\omega,z)|\leq |z|$, $\forall z\in\R$.

$\bullet$ \textit{Regularity of ${\bfm \gamma}$, $\nu$ and ${\bfm c}$}: Clearly, assumption 3) of Lemma \ref{constructiongamma} makes  Assumption \ref{regul}.2 hold.  Our purpose is now to check Assumptions \ref{regul}.3 and \ref{regul}.4.

For $|z|>0$ and each fixed $\omega$, the mapping $x\in\R\mapsto {\bfm h}(\tau_x\omega,z)$ is smooth (because of the regularity of ${\bfm c}^s$, see point 3) of Lemma \ref{constructiongamma}). From this and the relation ${\bfm h}(\omega, {\bfm \gamma}(\omega,z))=F(z)$, we let the reader deduce that the mapping $x\in\R\mapsto {\bfm \gamma}(\tau_x\omega,z) $ is also smooth.

By differentiating the relation $ {\bfm h}(\omega, {\bfm \gamma}(\omega,z))=F(z)$ with respect to $\omega$, we can compute the derivative $D{\bfm \gamma}$ 
$$D{\bfm \gamma}(\omega,z)=
\frac{ \int_{{\bfm \gamma}(\omega,z)}^{+\infty}D{\bfm c}^s(\omega,r)\chi(r)\,dr}{{\bfm c}^s(\omega,{\bfm \gamma}(\omega,z))\chi({\bfm \gamma}(\omega,z))}, \, \text{ if }z>0,\,\,\text{or }\frac{ \int_{{\bfm \gamma}(\omega,z)}^{+\infty}D{\bfm c}^s(\omega,r)\chi(r)\,dr}{{\bfm c}^s(\omega,{\bfm \gamma}(\omega,z))\chi({\bfm \gamma}(\omega,z))}, \, \text{ if }z>0. $$
For  $|{\bfm \gamma}(\omega,z)|\leq 1$, we can use point 1) of  Lemma \ref{constructiongamma}. Furthermore, we use the assumptions $|D{\bfm c}^s(\cdot,z)|_\infty\leq C_1$ and $0<m\leq {\bfm c}^s$ to deduce 
\begin{equation}\label{dgamma}
|D{\bfm \gamma}(\omega,z)|\one_{\{|{\bfm \gamma}(\omega,z)|\leq 1\}}\leq \frac{C_1M'}{m}|{\bfm \gamma}(\omega,z)|\one_{\{|{\bfm \gamma}(\omega,z)|\leq 1\}}.
\end{equation}
We are now in position to check Assumption \ref{regul}.3. By using the relation $|{\bfm \gamma}(\omega,z)|\leq |z|$ and \eqref{dgamma}, we have, for any $x,y\in\R$, 
\begin{align*}
\int_\R|{\bfm \gamma}&(\tau_x\omega,z)-{\bfm \gamma}(\tau_y\omega,z)|^2\one_{|z|\leq 1}\nu(z)\,dz\\
&\leq
\int_\R|y-x|^2\int_0^1|D{\bfm \gamma}(\tau_{(1-t)x+ty}\omega,z)|^2\,dt\one_{|z|\leq 1}\nu(z)\,dz\\
&\leq |y-x|^2\int_0^1\int_\R |D{\bfm \gamma}(\tau_{(1-t)x+ty}\omega,z)|^2\one_{\{|{\bfm \gamma}(\tau_{(1-t)x+ty}\omega,z)|\leq 1\}}\nu(z)\,dz\,dt\\
&\leq |y-x|^2 \frac{(C_1M')^2}{m^2}\int_0^1\int_\R|{\bfm \gamma}(\tau_{(1-t)x+ty}\omega,z)|^2\one_{|{\bfm \gamma}(\tau_{(1-t)x+ty}\omega,z)|\leq 1}\nu(z)dz\,dt\\ & \leq 
|y-x|^2 \frac{(C_1M')^2}{m^2}\int_0^1\int_\R |z|^2\one_{|z|\leq 1}{\bfm c}^s(\tau_{(1-t)x+ty}\omega,z)\chi(z)dz\,dt
\end{align*}
We easily conclude by using the bound $ {\bfm c}^s(\cdot,z)\leq M$ and $\int_\R \min(|z|^2,1) \chi(z)dz<+\infty$.
Finally, the relation $ |{\bfm \gamma}(\omega,z)|\leq |z|$ implies $\int_\R|{\bfm \gamma}(\tau_x\omega,z)|^2\one_{|z|\leq 1}\nu(z)\,dz\leq \int_\R z^2\one_{|z|\leq 1}\nu(z)\,dz$ so that we have checked Assumption \ref{regul}.3.

We now focus on Assumption \ref{regul}.4. First notice that the relation $|{\bfm \gamma}(\omega,z)|\leq |z|$ implies that the sets $\{z;|{\bfm \gamma}(\omega,z)|>1\}$ and $\{z;|z|\leq 1\}$ are disjoint. Hence, for $ \alpha>0$, we have
\begin{align}
\int_{\alpha<|{\bfm \gamma}(\omega,z)|}{\bfm \gamma}(\omega,z)\one_{|z|\leq 1}\nu(z)dz= &\int_{\alpha<|{\bfm \gamma}(\omega,z)|\leq 1}{\bfm \gamma}(\omega,z)\nu(z)dz- \int_{\alpha<|{\bfm \gamma}(\omega,z)|\leq 1}{\bfm \gamma}(\omega,z)\one_{|z|> 1}\nu(dz)\nonumber\\
=&\int_{\alpha<|z|\leq 1}z{\bfm c}^s(\omega,z)\chi(z)dz- \int_{\alpha<|{\bfm \gamma}(\omega,z)|\leq 1}{\bfm \gamma}(\omega,z)\one_{|z|> 1}\nu(z)dz\label{decomp}
\end{align}
Clearly, the second integral converges  towards $ \int_{|{\bfm \gamma}(\omega,z)|\leq 1}{\bfm \gamma}(\omega,z)\one_{|z|> 1}\nu(z)dz$ as $\alpha\to 0$ in $L^\infty(\Omega)$. Concerning the first integral, the convergence in $L^\infty(\Omega)$ is established in Lemma \ref{limitz} towards 
  $\frac{1}{2}\int_{|z|\leq 1}z^2\int_0^1D{\bfm c}^s(\tau_{rz}\omega,-z)\,dr\chi(z)dz $ as $\alpha\to 0$. Hence, we have proved that the following limit holds in $L^\infty(\Omega)$:
\begin{align*}
\lim_{\alpha\to 0} &\int_{\alpha<|{\bfm \gamma}(\omega,z)|}{\bfm \gamma}(\omega,z)\one_{|z|\leq 1}\nu(z)dz\\&=\frac{1}{2}\int_{|z|\leq 1}z^2\int_0^1D{\bfm c}^s(\tau_{rz}\omega,-z)\,dr\chi(z)dz- \int_{|{\bfm \gamma}(\omega,z)|\leq 1}{\bfm \gamma}(\omega,z)\one_{|z|> 1}\nu(z)dz.
\end{align*}
 It remains to prove that the limit satisfies a Lipschitz condition. From the regularity  of ${\bfm c}^s$, $\mu$ a.s., the mapping $x\in\R\mapsto \frac{1}{2}\int_{|z|\leq 1}z^2\int_0^1D{\bfm c}^s(\tau_{rz}\omega,-z)\,dr\chi(dz) $ is Lipschitzian. So, it just remains to prove that $\mu$ a.s., the mapping 
 $$\Gamma_\omega: x\in\R\mapsto  \int_{|{\bfm \gamma}(\tau_x\omega,z)|\leq 1}{\bfm \gamma}(\tau_x\omega,z)\one_{|z|> 1}\nu(dz)=\int_{z\in A(\tau_x\omega)}z{\bfm c}^s(\tau_x\omega,z)\chi(z)dz$$
  is Lipschitzian, where $$A(\omega)=\{z\in\R; |z|\leq 1\text{ and }z\not\in [{\bfm \gamma}(\omega,-1);{\bfm \gamma}(\omega,1)]\}.$$ For $x,y\in\R$, we define $A_{x,y}(\omega)$ as the symmetric difference of the sets $A(\tau_x\omega)$ and $A(\tau_y\omega)$: 
  $$A_{x,y}(\omega)=\big(A(\tau_x\omega)\setminus A(\tau_y\omega)\big)\cup \big(A(\tau_y\omega)\setminus A(\tau_x\omega)\big).$$
   For $z>0$, the relation $$F(z)={\bfm h}(\omega,{\bfm \gamma}(\omega,z))\geq \frac{m}{M}F({\bfm \gamma}(\omega,z))$$ leads to ${\bfm \gamma}(\omega,1)\geq F^{-1}(\frac{M}{m}F(1))$. Similarly, we have $ {\bfm \gamma}(\omega,-1)\leq F^{-1}(\frac{M}{m}F(-1))$. Hence, we can find $\beta>0$ such that $A(\omega)\subset \{z;\beta\leq |z|\leq 1\}$ for any $\omega\in \Omega$. Moreover, from \eqref{dgamma}, we have $|D{\bfm \gamma}(\omega,1)|\leq C_1M'/m$. In particular, the mapping $x\in\R\mapsto {\bfm \gamma}(\tau_x\omega,1)$ is $C_1M'/m$-Lipschitzian. It is plain to deduce that 
  $\int_{A_{x,y}(\omega)}dz\leq 2(C_1M'/m)|y-x|$. Finally, we conclude: for $x,y\in\R$, we have:
 \begin{align*}
|\Gamma_\omega(y)&-\Gamma_\omega(x)|\\
\leq &\int_{A(\tau_y\omega)}z|{\bfm c}^s(\tau_y\omega,z)-{\bfm c}^s(\tau_x\omega,z)|\chi(z)dz+\int_{\R}z{\bfm c}^s(\tau_x\omega,z)|\one_{A(\tau_y\omega)}-\one_{A(\tau_x\omega)}|\chi(z)dz\\
\leq & C_1|y-x|\int_{\beta\leq |z|\leq 1}\chi(z)dz+M\int_{\R}\one_{A_{x,y}(\omega)}\chi(z)dz\\
\leq & C_1|y-x|\int_{\beta\leq |z|\leq 1}\chi(z)dz+M\sup_{\beta\leq |z|\leq 1} \chi(z)2(C_1M'/m)|y-x|.
\end{align*}
Hence, the drift term $\lim_{\alpha\to 0} \int_{\alpha<|{\bfm \gamma}(\omega,z)|}{\bfm \gamma}(\omega,z)\one_{|z|\leq 1}\nu(z)dz$ is Lipschitzian.\qed

\vspace{2mm}
\noindent {\bf Proof of Lemma \ref{lemmapjs}.}
$\bullet$ \textit{Study of the convergence rate:} We have to compute the limit (in $L^1(\Omega)$) $$ \lim_{\epsilon\to 0}\frac{1}{\epsilon}\int_{\R}g(\delta(\epsilon)z){\bfm c}^s(\omega,z)\chi(dz) $$ for $g=\one_{[a,b]}$ such that 
$0\not \in[a,b]$. 

Since ${\bfm c}^s$ can be decomposed as ${\bfm c}^s(\omega,z)=\frac{1}{2}\big({\bfm c}(\tau_z\omega,-z)+{\bfm c}(\omega,z)\big)$, it suffices to compute the limits $ \lim_{\epsilon\to 0}\frac{1}{\epsilon}\int_{\R}g(\delta(\epsilon)z){\bfm c}(\omega,z)\chi(dz)$ and $ \lim_{\epsilon\to 0}\frac{1}{\epsilon}\int_{\R}g(\delta(\epsilon)z){\bfm c}(\tau_z\omega,-z)\chi(dz)$. The first limit raises no difficulty and matches $\int_{\R}g(z){\bfm \theta}(\omega,{\rm sign}(z))\mathcal{H}(dz)$ by using the convergence of ${\bfm c}$ (ass. 4 of Lemma \ref{lemmapjs}).

We now compute the second limit. By using the convergence of ${\bfm c}$ again and the invariance of the measure $\mu$ under $(T_z)_z$, one can establish
$$\lim_{\epsilon\to 0}\M\Big|\frac{1}{\epsilon}\int_{\R}g(\delta(\epsilon)z){\bfm c}(\tau_z\omega,-z)\chi(dz) -\frac{1}{\epsilon}\int_{\R}g(\delta(\epsilon)z){\bfm \theta}(\tau_z\omega,-{\rm sign}(z))\chi(dz)\Big|=0,$$
so that the proof boils down to establishing the following convergence
$$\lim_{\epsilon\to 0}\frac{1}{\epsilon}\int_{\R}g(\delta(\epsilon)z){\bfm \theta}(\tau_z\omega,-{\rm sign}(z))\chi(dz)=\int_{\R}g(z)\M[{\bfm \theta}(\omega,{\rm sign}(z))]\mathcal{H}(dz) .$$
Obviously, it suffices to establish that, for any function ${\bfm f}\in L^2(\Omega)$, 
\begin{align*}
\lim_{\epsilon\to 0}\frac{1}{\epsilon}\int_{\R}g(\delta(\epsilon)z){\bfm f}(\tau_z\omega)\chi(dz)=\int_{\R}g(z)\M[{\bfm f}]\mathcal{H}(dz),\quad \text{in }L^2(\Omega) 
\end{align*}
Actually this is a direct consequence of the spectral theorem. Let us explain why. Since $(T_z)_z$ is a strongly continuous group of unitary maps in $L^2(\Omega)$, there exists a projection valued measure  $E$ such that $(T_z{\bfm f},{\bfm g})_2=\int_\R e^{izu}E_{{\bfm f},{\bfm g}}(du)$, for any $z\in \R$ and ${\bfm f},{\bfm g}\in L^2(\Omega)$. Fix ${\bfm f}\in L^2(\Omega)$. Define the functions $b_\epsilon(u)= \frac{1}{\epsilon}\int_{\R}g(\delta(\epsilon)z)e^{izu}\chi(dz)$ $(\epsilon>0)$ and  the function $a(u)=\one_{u=0}\int_{\R}g(z)\mathcal{H}(dz)$ for $u\in\R$. Finally, set ${\bfm h}=\int_{\R}a(u)E_{{\bfm f}}(du)\in L^2(\Omega)$. Then
\begin{align*}
 \M\Big|\frac{1}{\epsilon}&\int_{\R}g(\delta(\epsilon)z){\bfm f}(\tau_z\omega)\chi(dz)-{\bfm h}\Big|^2\leq  \int_{\R}\big|b_\epsilon(u)-a(u)\big|^2E_{{\bfm f},{\bfm f}}(du)
 \end{align*}
 From the Lebesgue dominated convergence theorem, the last quantity tends to $0$ as $\epsilon\to 0$. Moreover, for any $z\in\R$, $T_z{\bfm h}=\int_{\R}e^{izu}a(u)E_{{\bfm f}}(du)=\int_{\R}a(u)E_{{\bfm f}}(du)={\bfm h}$, so that (by ergodicity of the measure $\mu$) ${\bfm h}=\M[{\bfm h}]=\M[{\bfm f}]\times\int_\R g(z)\mathcal{H}(dz)$.\qed

\section{Study of the Dirichlet form $B^s_\lambda $}\label{app:dirichlet}

This section is devoted to the  proofs of section \ref{sec:dirichlet}.

\vspace{2mm}
\noindent {\bf Proof of Lemma \ref{genito}.} Fix ${\bfm \varphi}\in C^2(\Omega)$.  The first step consists in computing $B^s({\bfm \varphi},{\bfm \psi})$ for any function ${\bfm \psi}\in \H$. To this purpose, first notice that an integration by parts yields:
$$B^d({\bfm \varphi},{\bfm \psi}) =\frac{1}{2}({\bfm a}D{\bfm \varphi},D{\bfm \psi})_\pi=\frac{1}{2}(e^{-2{\bfm V}}{\bfm a}D{\bfm \varphi},D{\bfm \psi})_2=-\frac{1}{2}(e^{2{\bfm V}}D(e^{-2{\bfm V}}{\bfm a}D{\bfm \varphi}),{\bfm \psi})_\pi.$$
Concerning $B^j$, by integrating by parts as in the proof of Lemma \ref{ippj}, we obtain:
$$B^j({\bfm \varphi},{\bfm \psi}) =-\lim_{\alpha\to 0}\int_{|z|\geq \alpha}\M\big[(T_z{\bfm \varphi}-{\bfm \varphi}){\bfm \psi}{\bfm c}(\cdot,z)\big]\chi(dz).$$ 
Notice that the existence of the limit raises no difficulty because of Lemma \ref{limitz} (take ${\bfm g}=(T_z{\bfm \varphi}-{\bfm \varphi}){\bfm \psi}$).
By using the relation $\nu\circ  {\bfm \gamma}_\omega^{-1}=e^{2{\bfm V}}{\bfm c}(\omega,z)\chi(dz)$, we deduce:
\begin{align*}
B^j({\bfm \varphi},{\bfm \psi}) 
=&-\lim_{\alpha\to 0}\M\big[\int_{|{\bfm \gamma}(\cdot,z)|\geq \alpha}(T_{{\bfm \gamma}(\cdot,z)}{\bfm \varphi}-{\bfm \varphi})\nu(dz){\bfm \psi}e^{-2{\bfm V}}\big]\\
=&-\M\big[\int_{\R}(T_{{\bfm \gamma}(\cdot,z)}{\bfm \varphi}-{\bfm \varphi}-{\bfm \gamma}(\cdot,z)\one_{|z|\leq 1}D{\bfm \varphi})\nu(dz){\bfm \psi}e^{-2{\bfm V}}\big]\\
&+\M\big[\lim_{\alpha\to 0}\int_{|{\bfm \gamma}(\cdot,z)|\geq \alpha}{\bfm \gamma}(\cdot,z)\one_{|z|\leq 1}\nu(dz)D{\bfm \varphi}{\bfm \psi}e^{-2{\bfm V}}\big]
\end{align*}
Gathering the above equalities, we have $B_\lambda^s({\bfm \varphi},{\bfm \psi}) =\lambda({\bfm \varphi},{\bfm \psi})_\pi-({\bfm L}'{\bfm \varphi},{\bfm \psi})_\pi$ (where ${\bfm L}'{\bfm \varphi}$ is given by the right-hand side of \eqref{generator}) for any function ${\bfm \psi}\in\H$. Hence, the mapping ${\bfm \psi}\in\H\mapsto B^s_\lambda({\bfm \varphi},{\bfm \psi})$ is $L^2(\Omega)$-continuous, ${\bfm \varphi}\in{\rm Dom}({\bfm L})$ and ${\bfm L}{\bfm \varphi}$ is given by \eqref{generator}.\qed

\vspace{3mm}
\noindent {\bf Proof of Proposition \ref{core}.} Let us first introduce the difference operator $\Gamma_r:L^p(\Omega)\to L^p(\Omega)$ $(r\not =0$ and $p\in[1,+\infty]$) defined by $\Gamma_r{\bfm \varphi}=\frac{1}{r}(T_r{\bfm \varphi}-{\bfm \varphi})$. It is straightforward to check the following properties:
\begin{align}\label{differ}
\forall {\bfm \varphi},{\bfm \psi}\in L^2(\Omega),\quad & (\Gamma_r{\bfm \varphi},{\bfm \psi})_2=-({\bfm \varphi},\Gamma_{-r}{\bfm \psi})_2\text{ and } \Gamma_r( {\bfm \varphi}{\bfm \psi})=T_r {\bfm \varphi}\Gamma_r{\bfm \psi}+ \Gamma_r{\bfm \varphi}{\bfm \psi}\\
\label{differ2} \forall{\bfm \varphi}\in{\rm Dom}(D),\quad & |\Gamma_r{\bfm \varphi}|_p\leq |D{\bfm \varphi}|_p.
\end{align}
Fix ${\bfm f}\in L^2(\Omega)$ and denote $G_\lambda{\bfm f}$ by ${\bfm f}_\lambda$. 
%

Choose ${\bfm \psi}\in\H$ and $r\not =0$, an plug $\Gamma_r{\bfm \psi}$ into \eqref{eqbase}:
\begin{align*}
({\bfm f},\Gamma_r{\bfm \psi})_\pi=&\lambda({\bfm f}_\lambda,\Gamma_r{\bfm \psi})_\pi+\frac{1}{2}({\bfm a}D{\bfm f}_\lambda,D\Gamma_r{\bfm \psi})_\pi+\frac{1}{2}\M\int_{\R}(T_z{\bfm f}_\lambda-{\bfm f}_\lambda)(T_z\Gamma_r{\bfm \psi}-\Gamma_r{\bfm \psi}){\bfm c}(\cdot,z)\chi(dz)
\end{align*}
Then we use \eqref{differ} to obtain
\begin{align}\label{eqr}
({\bfm f},\Gamma_r{\bfm \psi})_\pi=&-B_\lambda^s(\Gamma_{-r}{\bfm f}_\lambda,{\bfm \psi})-\lambda(\Gamma_{-r}(e^{-2{\bfm V}})T_{-r}{\bfm f}_\lambda,{\bfm \psi})_2-\frac{1}{2}(\Gamma_{-r}({\bfm a}e^{-2{\bfm V}})DT_{-r}{\bfm f}_\lambda,D{\bfm \psi})_2\nonumber\\
&-\frac{1}{2}\M\int_{\R^*}(T_zT_{-r}{\bfm f}_\lambda-T_{-r}{\bfm f}_\lambda)(T_z{\bfm \psi}-{\bfm \psi})\Gamma_{-r}{\bfm c}(\cdot,z)\chi(dz),
\end{align}
From  estimate \eqref{differ2} and Lemma \ref{contdc} (take ${\bfm d}=\Gamma_{-r}{\bfm c}$ and $M=C_1$, $C_1$ given by Assumption \ref{regul}), we deduce 
\begin{align*}B^s_\lambda(\Gamma_{-r}{\bfm f}_\lambda,{\bfm \psi})\leq &|e^{-2{\bfm V}}|_\infty|{\bfm f}|_2|D{\bfm \psi}|_2+\lambda|D(e^{-2{\bfm V}})|_\infty|{\bfm f}_\lambda|_2|{\bfm \psi}|_2+\frac{1}{2}|D({\bfm a}e^{-2{\bfm V}})|_\infty|D{\bfm f}_\lambda|_2\\&+C_{\ref{contdc}}|D{\bfm f}_\lambda|_2||D{\bfm \psi}|_2+C_{\ref{contdc}}|{\bfm f}_\lambda|_2||{\bfm \psi}|_2\\
\leq & C(|{\bfm \psi}|_\pi+|D{\bfm \psi}|_\pi)
\end{align*}
where the constant $C$ does not depend on $r$ (only on the regularity of ${\bfm a}$, ${\bfm V}$, ${\bfm c}$, on $\chi $ and on the norms $|{\bfm f}_\lambda|_2 $ and $|D{\bfm f}_\lambda|_2 $).
Choosing ${\bfm \psi}=\Gamma_{-r}{\bfm f}_\lambda$ in the previous inequality yields
\begin{align*}
B^s_\lambda(\Gamma_{-r}{\bfm f}_\lambda,\Gamma_{-r}{\bfm f}_\lambda)\leq C(|\Gamma_{-r}{\bfm f}_\lambda|_\pi+|D\Gamma_{-r}{\bfm f}_\lambda|_\pi)\leq \frac{C^2}{2\lambda}+\frac{\lambda}{2}|\Gamma_{-r}{\bfm f}_\lambda|_\pi^2+C^2M_{\ref{ellipticity}}+\frac{M^{-1}_{\ref{ellipticity}}}{4}|D\Gamma_{-r}{\bfm f}_\lambda|_\pi^2,
\end{align*}
 in such a way that $B^s_\lambda(\Gamma_{-r}{\bfm f}_\lambda,\Gamma_{-r}{\bfm f}_\lambda)\leq  \frac{C^2}{\lambda}+2C^2M_{\ref{ellipticity}}$. Hence, the family $(\Gamma_{-r}{\bfm f}_\lambda)_{r\not =0}$ is bounded in $\H$, and is therefore  weakly compact in $\H$. By passing to the limit in \eqref{eqr} as $r\to 0$, it is plain to see that the limit  ${\bfm g}_\lambda \in\H$ (in fact ${\bfm g}_\lambda=D{\bfm f}_\lambda $) of a converging subsequence satisfies the relation (for each $ {\bfm \psi}\in \H$)
 \begin{align}\label{eqli}
({\bfm f},D{\bfm \psi})_\pi=&-B_\lambda^s({\bfm g}_\lambda,{\bfm \psi})-\lambda(D(e^{-2{\bfm V}}){\bfm f}_\lambda,{\bfm \psi})_2-\frac{1}{2}(D({\bfm a}e^{-2{\bfm V}})D{\bfm f}_\lambda,D{\bfm \psi})_2\nonumber\\
&-\frac{1}{2}\M\int_{\R}(T_z{\bfm f}_\lambda-{\bfm f}_\lambda)(T_z{\bfm \psi}-{\bfm \psi})D{\bfm c}(\cdot,z)\chi(dz).
\end{align}
(The $\H$-continuity of the last integral is proved in Lemma \ref{contdc} with ${\bfm d}=D{\bfm c}$.) In particular, the relation $D{\bfm f}_\lambda={\bfm g}_\lambda\in \H$  implies  that $D{\bfm f}_\lambda\in{\rm Dom}(D)$.  We have proved $G_\lambda(L^2(\Omega))\subset H^2(\Omega)$.

We prove now that ${\bfm f}\in H^1(\Omega)\Rightarrow {\bfm f}_\lambda\in H^3(\Omega)$. If we can prove that $D{\bfm f}_\lambda$ is the solution to an equation of the type $B_\lambda^s(D{\bfm f}_\lambda,{\bfm \psi})=({\bfm g},{\bfm \psi})_\pi $ with ${\bfm g}\in L^2(\Omega)$, then $ D{\bfm f}_\lambda\in H^2(\Omega)$ (i.e. ${\bfm f}_\lambda\in H^3(\Omega)$) according to the previous argument. That is what we are going to prove.  In the case ${\bfm f}\in H^1(\Omega) $, equation \eqref{eqli} becomes (by integrating by parts in \eqref{eqli} the terms containing $D{\bfm \psi}$ and by using Lemma \ref{contbr} below)
\begin{align}\label{eqli2}
B_\lambda^s(D{\bfm f}_\lambda,{\bfm \psi})=&(e^{2{\bfm V}}D(e^{-2{\bfm V}}{\bfm f}),{\bfm \psi})_\pi-\lambda(e^{2{\bfm V}}D(e^{-2{\bfm V}}){\bfm f}_\lambda,{\bfm \psi})_\pi+\frac{1}{2}\big(e^{2{\bfm V}}D(D({\bfm a}e^{-2{\bfm V}})D{\bfm f}_\lambda),{\bfm \psi}\big)_\pi\nonumber\\
&+\M\big[\lim_{\alpha\downarrow 0}\int_{|z|>\alpha}(T_z{\bfm f}_\lambda-{\bfm f}_\lambda)e^{2{\bfm V}}D{\bfm c}(\cdot,z) \chi(dz){\bfm \psi}e^{-2{\bfm V}}\big].
\end{align}
So we have ${\bfm g}=e^{2{\bfm V}}\Big(D(e^{-2{\bfm V}}{\bfm f}) -\lambda D(e^{-2{\bfm V}}){\bfm f}_\lambda+\frac{1}{2}D(D({\bfm a}e^{-2{\bfm V}})D{\bfm f}_\lambda)+\lim_{\alpha\downarrow 0}\int_{|z|\geq \alpha}(T_z{\bfm f}_\lambda-{\bfm f}_\lambda)D{\bfm c}(\cdot,z) \chi(dz)\Big)$. 

As guessed by the reader, the proof is now completed recursively, the only difficulty being of notational nature. \qed


\begin{lemma}\label{contbr}
For any ${\bfm f}\in H^2(\Omega)$ and ${\bfm g}\in\H$, the following integration by parts holds:
$$-\frac{1}{2}\M\int_{\R}(T_z{\bfm f}-{\bfm f})(T_z{\bfm g}-{\bfm g})D{\bfm c}(\cdot,z)\chi(dz)=\M\big[\lim_{\alpha\downarrow 0}\int_{|z|>\alpha}(T_z{\bfm f}-{\bfm f})e^{2{\bfm V}}D{\bfm c}(\cdot,z) \chi(dz){\bfm g}e^{-2{\bfm V}}\big].$$
\end{lemma}

\noindent {\bf Proof.} First notice that the limit $\lim_{\alpha\downarrow 0}\int_{|z|>\alpha}(T_z{\bfm f}-{\bfm f})e^{2{\bfm V}}D{\bfm c}(\cdot,z) \chi(dz)$ is well defined in the $L^2(\Omega)$ sense thanks to Lemma \ref{limitz} (take ${\bfm g}=T_z{\bfm f}-{\bfm f}$). To prove the integration by parts formula above, we can make the same computations as in the proof of Lemma \ref{ippj} (use the symmetry of $D{\bfm c}$). Details are left to the reader. We also point out that the same property holds for the successive derivatives of ${\bfm c}$.\qed

\vspace{2mm}
\noindent {\bf Proof of Proposition \ref{submarkov}.} The proof is not specifically written for a random medium. However, the arguments used in \cite{fuku} do not fail in our framework. It suffices to prove that the symmetric form $B^d$ is Markovian (cf \cite{fuku}), which can be established by following the proofs of \cite{fuku} or \cite[examples 3.6.8 and 3.6.9]{applebaum}.\qed 

\end{document}